\documentclass[12pt]{article}

\usepackage{amsmath,enumerate,amssymb,amsthm,color}
\usepackage{pgf,tikz}
\usepackage{color}

\title{Isomorphic tetravalent cyclic Haar graphs} 
\author{ S. Hiroki  Koike-Quintanar\footnote{Supported in part by ARRS -- 
Agencija za raziskovanje Republike Slovenije, program no. P1-0285. } \\  
{\small University of Primorska, UP IAM} \\  [-0.8ex]
{\small  Muzejski trg 2, SI6000 Koper, Slovenia} \\ 
{\small \texttt{hiroki.koike@upr.si}} 
\and    
Istv\'an Kov\'acs\footnote{Supported in part by ARRS -- 
Agencija za raziskovanje Republike Slovenije, program no. P1-0285.} \\   
{\small University of Primorska, UP IAM and UP FAMNIT}  \\  [-0.8ex]
{\small  Muzejski trg 2, SI6000 Koper, Slovenia} \\ 
{\small  \texttt{istvan.kovacs@upr.si}}
}

\newtheorem{thm}{Theorem}[section]
\newtheorem{lem}[thm]{Lemma}
\newtheorem{cor}[thm]{Corollary}
\newtheorem{prop}[thm]{Proposition}
\theoremstyle{definition}
\newtheorem{exm}[thm]{Example}

\def\conga{\cong_{\mathrm{aff}}}
\def\congc{\cong_{\mathrm{cay}}}
\def\I{\mathcal{C}_{\mathrm{iso}}}
\def\AI{\mathcal{C}_{\mathrm{aff}}}
\def\Z{\mathbb{Z}}
\def\GBG{\mathbb{B}} 
\def\proof{\noindent{\sc Proof}.\ }
\def\QED{\hfill\vrule height .9ex width .8ex depth -.1ex \medskip}
\DeclareMathOperator{\aut}{Aut}
\DeclareMathOperator{\cay}{Cay}
\DeclareMathOperator{\iso}{Iso}
\DeclareMathOperator{\orb}{Orb}
\DeclareMathOperator{\sym}{Sym}
\newcommand{\sg}[1]{\langle {#1}\rangle}
\newcommand{\comment}[1]{}

\begin{document}

\maketitle

\begin{abstract}
Let  $S$ be  a  subset  of the cyclic group $\Z_n$. The cyclic Haar graph $H(\Z_n,S)$ is the bipartite graph with color classes 
$\Z_n^+$ and $\Z_n^-,$ and edges $\{x^+,y^-\},$ where $x,y \in \Z_n$ and $y - x \in S$.  
In this paper we give sufficient and necessary conditions for the isomorphism of two connected cyclic Haar graphs 
of valency $4$.  


\end{abstract}

\section{Introduction}

Let $S$ be a subset of a finite group $G$. The {\em Haar graph}  
$H(G,S)$ is the bipartite graph with color classes  identified with $G$ and written as 
$G^+$ and $G^-,$ and the edges are $\{x^+,y^-\},$ where $x,y \in G$ and $y x^{-1}  \in S$.  
Haar graphs were introduced for abelian groups by Hladnik, Maru\v si\v c and Pisanski \cite{HlaMP02}, and 
were redefined under the name {\em bi-Cayley graphs} in  \cite{XuJSL08}.  A Haar graph $H(G,S)$  is called {\em cyclic} if 
$G$ is a cyclic group. In this paper we consider the problem of giving sufficient and necessary conditions for the 
isomorphism of two cyclic Haar graphs. This is a natural continuation of 
the  isomorphism problem of circulant digraphs which has been solved by Muzychuk \cite{Muz04}.  It appears 
in the context of circulant matrices under the name 
bipartite \'Ad\'am problem \cite{WieZ07}, and  also  in the context of cyclic configurations \cite{BobPZ05,HlaMP02}.

The symbol $\Z_n$ denotes the additive group of the ring $\Z/n \Z$ of residue classes modulo $n,$ and 
$\Z_n^*$ denotes the multiplicative group of units in $\Z/n \Z$. Two Haar graphs 
$H(\Z_n,S)$ and $H(\Z_n,T)$ are called {\em affinely equivalent,} written as $H(\Z_n,S) \conga H(\Z_n,T),$ 
if $S$ can be mapped to $T$ by an affine transformation, i.e., $a S + b = T$ for some $a \in \Z_n^*$ and $b \in \Z_n$. 
It is an easy exercise to  show that two affinely equivalent cyclic Haar graphs are 
isomorphic as usual graphs. The converse implication is not true in general, and this makes the following definition interesting (see \cite{XuJSL08}): we say that a subset $S \subseteq \Z_n$
is a {\em BCI-subset} if for each $T \subseteq G,$ 
$H(\Z_n,S) \cong H(\Z_n,T)$  if and only if $H(\Z_n,S) \conga H(\Z_n,T)$.
Wiedemann and Zieve proved in \cite[Theorem 1.1]{WieZ07} that any subset $S$ of $\Z_n$ is a BCI-subset if $|S| \le 3$ 
(a special case was proved earlier  in \cite{BorBF73}).  However, this is not true if $|S| \ge 4$ (see  \cite{HlaMP02,WieZ07}), 
hence the first nontrivial case of the isomorphism problem occurs when $|S|=4$. In this paper we settle this case by proving the following theorem:

\begin{thm}\label{MAIN1}
Two connected Haar graphs $H(\Z_n,S)$ and $H(\Z_n,T)$ with $|S| = |T| =4$ are isomorphic if and only if there exist 
$a_1,a_2 \in \Z_n^*$  and $b_1,b_2 \in \Z_n$ such that 
\begin{enumerate}[(1)]
\item $a_1 S + b_1 = T;$  or 
\item $a_1 S + b_1 = \{ 0,u,v,v+m \}$  and $a_2 T + b_2 = \{ 0,u+m,v,v+m\},$ where $n=2m,$, $\Z_n = \langle u,v \rangle$, 
$2 \mid u,$ $2u \mid m$ and $u/2 \not\equiv v + m/(2u) (\mathrm{mod\;} m/u)$.
\end{enumerate}
\end{thm}

\noindent {\bf Remark.}   
A group $G$ is called an {\em $m$-BCI-group} if every subset $S$ of $G$ with $|S| \le m$ is  a BCI-subset (see \cite{JinL10, XuJSL08}). In this context \cite[Theorem 1.1]{WieZ07}  can be rephrased as $\Z_n$ is a $3$-BCI-group 
for any number $n$; and Theorem 1.1 says that $\Z_n$ is not a $4$-BCI-group if and only if  $n$ is divisible by $8$. 
This refines \cite[Theorem 7.2]{WieZ07} in which it is proved  that, if $\Z_n$ contains a  non-BCI-subset of size $k, k \in \{4,5\},$ 
then $n$ has a prime divisor  less or equal to  $2k(k-1)$. 

\medskip

Our approach is group theoretical, we adopt the ideas of  \cite{Bab77,Muz99}.
In short terms the initial problem is transformed to a problem about the automorphism group of the graphs in question. 
Theorem 1.1 is proven in two steps: first it is settled for graphs $H(\Z_n,S)$ with $S$  satisfying additional conditions (see Theorem 3.1); then  it is  shown that, if $S$ is not a BCI-subset, then it is affinely equivalent to a set satisfying the conditions of
 Theorem 3.1 (see Theorem 4.1).

\section{A Babai type theorem}

In this paper every group, graph and digraph is finite. For a (di)graph  $\Gamma,$ the symbols $V(\Gamma),$ $E(\Gamma)$ and $\aut(\Gamma)$ denote the set of its vertices, (directed) edges and the full group of  its automorphisms, respectively. 
Regarding terminology and notation in permutation group theory we follow \cite{DixM96}.  

Let $S$ be a subset of a group $G$. The {\em Cayley digraph} $\cay(G,S)$ is the digraph with 
vertex set $G,$ and its directed edges are $(x,y),$ where $x,y \in G$ and $y x^{-1} \in S$. 
Two digraphs $\cay(G,S)$ and $\cay(G,T)$ are called {\em Cayley isomorphic}, written as 
$\cay(G,S) \congc \cay(G,T)$, if $T = S^\varphi$ for  
some group automorphism $\varphi \in \aut(G)$. It is clear that such an automorphism induces an isomorphism between 
$\cay(G,S)$ and $\cay(G,T),$ and thus Cayley isomorphic digraphs are isomorphic as usual digraphs. It is also well-known that
the converse implication is not true, and this makes sense for the following 
definition (see \cite{Bab77}): a  subset $S \subseteq G$ is a {\em CI-subset} if for each $T \subseteq G,$ 
$\cay(G,S) \cong \cay(G,T)$  if and only if $\cay(G,S) \congc \cay(G,T)$.
The following equivalence was proved by Babai \cite[3.1 Lemma]{Bab77}.

\begin{thm}\label{B}
The following are equivalent for every Cayley digraph $\cay(G,S)$.
\begin{enumerate}[(1)]
\item $S$ is a CI-subset.
\item Every two regular subgroups of $\aut(\cay(G,S))$ isomorphic to $G$ are 
conjugate  in $\aut(\cay(G,S))$.
\end{enumerate}
\end{thm}

\noindent Theorem \ref{B} essentially says that the CI-property of a given subset $S$ depends entirely on the automorphism 
group $\aut(\cay(G,S))$.  In this section we  prove analogous results for cyclic Haar graphs.

Let $V=V(H(\Z_n,S))$ be the vertex set of the Haar Graph $H(\Z_n,S)$. Throughout this paper $c$ and $d$ denote the permutations of V defined by 
\begin{equation}\label{c d}
c : x^\varepsilon  \mapsto (x+1)^\varepsilon \text{ and } d :  x^\varepsilon \mapsto \begin{cases} (n-x)^-  & \text{ if } 
\varepsilon = +,  \\ (n-x)^+ & \text{ if } \varepsilon = -,\end{cases}
\end{equation}
where $x \in \Z_n$ and $\varepsilon \in \{+,-\}$.
It follows immediately that both $c$ and $d$ are automorphisms of any Haar graph $H(\Z_n,S)$. 
Denote by $C$ the group generated by $c,$ and by $D$ the group generated by $c$ and $d$. 
The group $D$ acts regularly on $V,$ and $D$  is isomorphic to $D_{2n}$. 
Thus $H(\Z_n,S)$ is isomorphic to a Cayley graph over $D,$ and so Theorem \ref{B} can be applied. 
The following  corollary is obtained.

\begin{cor}\label{corB}
The implication (1) $\Leftarrow$ (2) holds for every Haar graph $H(\Z_n,$ $S)$.  
\begin{enumerate}[(1)]
\item $S$ is a BCI-subset.  
\item Every two regular subgroups of $\aut(H(\Z_n,S))$ isomorphic to $D$ are conjugate in $\aut($ $H(\Z_n,S))$. 
\end{enumerate}
\end{cor}

However, we do not have equivalence in Corollary \ref{corB}  as it is shown in the following example. 

\begin{exm}\label{Haar10} 
Let $\Gamma = H(\Z_{10},\{0,1,3,4\})$. Using the computer package {\sc Magma} \cite{BosCP97} we 
compute that $\Gamma$ is edge-transitive and  its automorphism group $\aut(\Gamma) \cong D_{20} \rtimes \Z_4$. 
Furthermore, $\aut(\Gamma)$ contains a regular subgroup $X$ which is isomorphic to $D_{20}$ but $X \ne D_{20}$, 
hence (2) in Corollary 2.2 does not hold. 

On the other hand, we find that for every subset $T \subseteq \Z_{10}$ with $0 \in T$ and $|T| =4,$  
the corresponding Haar graph $H(\Z_{10},T) \cong \Gamma$ exactly when $H(\Z_{10},T) \conga \Gamma$. 
Thus $\{0,1,3,4\}$ is a BCI-subset, so (1) in Corollary 2.2 is true.  
\QED 
\end{exm}

Example \ref{Haar10} shows that the isomorphism problem of cyclic Haar graphs is not a particular case of the isomorphism problem of Cayley graphs over dihedral groups. We remark that the latter problem is still unsolved, for partial solutions, see  \cite{Bab77,QuY97,Spi08}. Nonetheless, the idea of Babai works well if  instead of the regular subgroup $D$ we consider its index $2$ cyclic subgroup $C$.

We say that  a  permutation group $G \le \sym(\Z_n^+ \cup \Z_n^-)$  is {\em bicyclic} if  
$G$ is a cyclic group which has two orbits: $\Z_n^+$ and $\Z_n^-$.  By a bicyclic group of a Haar graph $\Gamma = H(\Z_n,S)$  we simply mean a bicyclic subgroup $X \le \aut(\Gamma)$.
Obviously, $C$ is a bicyclic group of any cyclic Haar graph, and being so it will be referred to 
as the {\em canonical bicyclic group}.  

Let $\iso(\Gamma)$ denote the set of all isomorphisms from $\Gamma$ to any other 
Haar graph $H(\Z_n,T),$ i.e.,   
\[ \iso(\Gamma) = \big\{ f \in \sym(V) :  \Gamma^f = H(\Z_n,T) \text{ for some }  T \subseteq \Z_n \big\}.\] 
And let $\I(\Gamma)$ denote the {\em isomorphism class} of cyclic Haar graphs which contains $\Gamma,$ i.e.,  
$ \I(\Gamma) = \{ \Gamma^f : f \in \iso(\Gamma) \}$. 

\begin{lem}\label{Lem1}
Let $\Gamma = H(\Z_n,S)$  be a  connected  Haar graph and $f$ be in $\sym(V)$. Then 
$f \in \iso(\Gamma)$ if and only if $f  C f^{-1}$ is a bicyclic group of $\Gamma$.
\end{lem}

\proof 
Let $f \in \iso(\Gamma)$. Then $f C f^{-1} \le \aut(\Gamma)$. Clearly, $f C f^{-1}$ is a cyclic group. 
Since the sets $\Z_n^+$ and $\Z_n^-$ are the color classes of the connected bipartite graph $\Gamma,$ 
$f$ preserves these color classes, implying that 
$\orb(f C f^{-1},V) = \{ Z_n^+, \Z_n^- \}$. The group $f C f^{-1}$ is a bicyclic group of $\Gamma$.
 
Conversely, suppose that $f C f^{-1}$ is a bicyclic group of $\Gamma$. 
Then  $C =  f^{-1}(f C f^{-1})f  \le \aut( \Gamma^f)$. 
Because that $\orb(f C f^{-1},V) = \{ \Z_n^+, \Z_n^-\}$, the graph $\Gamma^f$ is  connected  and bipartite 
with color classes $\Z_n^+$ and $Z_n^-$. 
We conclude that $\Gamma^f = H(\Z_n,T)$ for some $T \subseteq \Z_n,$  so $f \in \iso(\Gamma)$. 
The lemma is proved. 
\QED

Lemma \ref{Lem1} shows that the normalizer $N_{\sym(V)}(C) \subseteq \iso(H(\Z_n,S))$. 
The group $N_{\sym(V)}(C)$ is known to consist of the  following permutations:  
\begin{equation}\label{varphi psi}
\varphi_{r,s,t} : x^\varepsilon \mapsto 
\begin{cases} 
(r x+s)^+ & \text{ if } \varepsilon = +, \\
(r x+t)^- & \text{ if } \varepsilon = -,
\end{cases}   \quad 
\psi_{r,s,t} : x^\varepsilon  \mapsto 
\begin{cases} 
(r x+s)^- & \text{ if } \varepsilon = +, \\
(r x+t)^+ & \text{ if } \varepsilon = -,
\end{cases}
\end{equation}
where $r \in \Z_n^*$ and $s,t \in \Z_n$. 
Note that, two Haar graphs $H(\Z_n,S)$ and $H(\Z_n,T)$ are from the same orbit under $N_{\sym(V)}(C)$ exactly when 
$H(\Z_n,S) \conga H(\Z_n,T)$.   
Let $\AI(\Gamma)$ denote the {\em affine equivalence class} of cyclic Haar graphs which contains the graph 
$\Gamma = H(\Z_n,S)$, i.e., $\AI(\Gamma) = \{  \Gamma^\varphi : \varphi \in N_{\sym(V)}(C) \}$. 
It is clear that the  isomorphism class $\I(\Gamma)$ splits into affine equivalence classes:
\[ \I(\Gamma) = \AI(\Gamma_1) \, \dot{\cup} \, \cdots \, \dot{\cup} \, \AI(\Gamma_k)
\footnote{Here we mean that 
$\I(\Gamma) = \AI(\Gamma_1) \, \cup \, \cdots \, \cup \, \AI(\Gamma_k)$ and 
$\AI(\Gamma_i) \cap \AI(\Gamma_j) = \emptyset$ for every $i,j \in \{1,\dots,k\},$ $i \ne j$.}.
\] 

Our next goal is to describe the above decomposition with the aid of bicyclic groups. 
Let  $X$ be a bicyclic group of a connected graph $\Gamma = H(\Z_n,S)$. Then $g^{-1} X g$ is also a 
bicyclic group for every $g \in \aut(\Gamma)$, hence the full set of bicyclic groups of $\Gamma$ is the union of $\aut(\Gamma)$-conjugacy classes.  We say that a subset $\Xi \subseteq \iso(\Gamma)$ is a {\em bicyclic base} of $\Gamma$ if 
the subgroups $\xi C \xi^{-1}, \xi \in \Xi,$ form a complete set of representatives of the corresponding conjugacy classes. 
Thus every bicyclic group $X$ can be expressed as    
\[X  = g \xi C (g \xi)^{-1}  \text{ for a unique } \xi \in \Xi \text{ and } g  \in \aut(\Gamma).\]

\noindent {\bf Remark.} 
Our definition of a bicyclic base copies in a sense the definition of a cyclic base introduced by 
Muzychuk  \cite[Definition, page 591]{Muz99}. 

\begin{thm}\label{Thm1}
Let $\Gamma = H(\Z_n,S)$  be a  connected  Haar graph with a bicyclic base $\Xi$. 
Then $\I(\Gamma) = \dot{\bigcup}_{\xi \in \Xi} \AI ( \Gamma^\xi )$.
\end{thm}

\proof  It follows immediately that,
\begin{equation}\label{subset} 
\I(\Gamma) \supseteq \bigcup_{\xi \in \Xi} \AI( \Gamma^\xi ).
\end{equation} 
We prove that equality holds in \eqref{subset}. Pick $\Sigma \in \I(\Gamma)$. Then $\Sigma = \Gamma^f$ for some 
$f \in \iso(\Gamma)$. By Lemma \ref{Lem1}, $f C f^{-1}$  is a bicyclic group of $\Gamma,$ hence 
 \[ f C f^{-1} = g \xi C (g \xi)^{-1}, \; \xi \in \Xi, g \in \aut(\Gamma). \]
Thus $ f^{-1} g \xi = h,$ where $h \in N_{\sym(V)}(C)$. Then 
\[ \Sigma = \Gamma^f = \Gamma^{g \xi h^{-1}} = \big( \Gamma^\xi \big)^{h^{-1}}. \]
This shows that $\Sigma \in \AI( \Gamma^\xi),$ and so   
\[ \I(\Gamma) \subseteq \bigcup_{\xi \in \Xi} \AI( \Gamma^\xi ). \]
In view of \eqref{subset} the  two sides are equal.  

Moreover, if $\Sigma \in \AI(\Gamma^{\xi_1}) \cap \AI(\Gamma^{\xi_2})$ for $\xi_1,\xi_2 \in \Xi,$ then  
$\Gamma^{\xi_1} = \Gamma^{\xi_2 h}$ for some $h \in N_{\sym(V)}(C)$. Hence  
$\xi_2 h \xi_1^{-1} = g$ for some $g \in \aut(\Gamma),$ and so 
\[
\xi_1 C  \xi_1^{-1} = g^{-1 }\xi_2 h C  h^{-1 }\xi_2^{-1} g = g^{-1 } ( \xi_2 C \xi_2^{-1}) g.
\]
The bicyclic subgroups $\xi_1 C  \xi_1^{-1}$ and $\xi_2 C  \xi_2^{-1}$ are conjugate in $\aut(\Gamma),$ hence 
$\xi_1 = \xi_2$ follows from the definition of the bicyclic base $\Xi$. 
We obtain that $\AI(\Gamma^{\xi_1}) \cap \AI(\Gamma^{\xi_2}) = \emptyset$ whenever 
$\xi_1,\xi_2 \in \Xi,$ $\xi_1 \ne \xi_2,$ and so $\I(\Gamma) = \dot{\bigcup}_{\xi \in \Xi} \AI ( \Gamma^\xi )$. The theorem is proved. 
\QED

As a direct consequence  of Theorem \ref{Thm1} we obtain the following corollary, analog of Theorem \ref{B}.   

\begin{cor}\label{corBCI}
The following are equivalent for every connected Haar graph $H(\Z_n,S)$.
\begin{enumerate}[(1)]
\item $S$ is a BCI-subset.
\item Any two bicyclic groups of $H(\Z_n,S)$ are conjugate in $\aut(H(\Z_n,S))$. 
\end{enumerate}
\end{cor}

In our last proposition we connect the BCI-property with the CI-property.    
For $a^\varepsilon \in V,$ in what follows $\aut(H(\Z_n,S))_{a^\varepsilon}$ denotes the vertex stabilizer of 
$a^\varepsilon$ in $\aut(H(\Z_n,S))$. 

\begin{prop}\label{Prop1}
Suppose that $\Gamma = H(\Z_n,S)$  is a connected Haar graph such that 
$\aut(\Gamma)_{0^+} = \aut(\Gamma)_{a^-}$ for some $a \in \Z_n$. Then the following are equivalent.
\begin{enumerate}[(1)]
\item $S$ is a BCI-subset.
\item $S-a = \{s-a : s \in S\}$ is a CI-subset. 
\end{enumerate}
\end{prop}

\proof For sake of simplicity we put  $A = \aut(\Gamma)$ and $G$ for the setwise stabilizer 
$\aut(\Gamma)_{\{\Z_n^+\}}$ of the color class $\Z_n^+$ in $\aut(\Gamma)$. Obviously, $X \le G$ for 
every bicyclic group $X$ of $\Gamma$. Since $A = G \rtimes \sg{d}$ and 
$d$ normalizes $C$, it follows that the conjugacy class of subgroups of $A$ containing $C$ is equal to 
the conjugacy class of subgroups of $G$ containing $C$. 
Using this and Theorem \ref{Thm1}, we obtain that $S$ is a BCI-subset if and only if   
every bicyclic group is conjugate to $C$ in $G$.  

Let $W = \{0^+,a^-\}$ and consider the setwise stabilizer $A_{\{W\}}$. 
Since $A_{0^+} = A_{a^-},$ $A_{0^+} \le A_{\{W\}}$.
By \cite[Theorem 1.5A]{DixM96},  the orbit of $0^+$ under  $A_{\{W\}}$ is a block of imprimitivity (for short a block) for $A$. 
Denote this block by $\Delta$ and the induced system of blocks by $\delta$ (i.e., $\delta = \{ \Delta^g : g \in G\}$). 
Consider the element $g = dc^a$ from $D$. We see that $g$ switches $0^+$ and $a^-,$ hence $A_{\{W\}} = A_{0^+} 
\big\langle g \big\rangle$. Therefore, $\Delta =  (0^+)^{A_{\{W\}}} = (0^+)^{A_{0^+} \langle g \rangle} = (0^+)^{\langle g \rangle} = W,$ and so  
\[  \delta = \big\{ \, \{x^+, (x+a)^-\} : x \in \Z_n \, \big\}. \]

Define the mapping 
$\varphi : \delta \to \Z_n$  by $\varphi : \{x^+, (x+a)^-\} \mapsto x,$ $x \in \Z_n$. Now, an 
action of $A$ on $\Z_n$ can be defined by letting $g \in A$ act as 
\[ x^g = x^{\varphi^{-1} g \varphi}, \; x \in \Z_n. \]
For $g \in A$ we write $\bar{g}$ for the image of $g$ under the corresponding  permutation 
representation, and for a subgroup $X \le A$ we let $\bar{X} = \{\bar{x} : x \in X\}$. In this action of $A$ the 
subgroup $G < A$ is faithful. Also notice that,  a subgroup $X \le G$ is a bicyclic group of $\Gamma$ 
if and only if $\bar{X}$ is a regular cyclic subgroup of $\bar{G}$. 
In particular, for the canonical bicyclic group $C,$
$\bar{C} = (\Z_n)_{\rm right},$ where $(\Z_n)_{\rm right}$ denotes the group generated by 
the affine transformation $x \mapsto x+1,$ $x \in \Z_n$.

Pick $g \in G$ and $(x,x+s-a) \in \Z_n \times \Z_n,$ where $s \in S$.
Then $g$ maps the directed edge $(x^+,(x+s)^-)$ to a directed edge $(y^+,(y+q)^-)$ for some $y \in \Z_n$ and 
$q\in S$. Since $\delta$ is a system of blocks for $G,$ $g$ maps $(x+s-a)^+$ to $(y+q-a)^+$, and so $\bar{g}$ maps the pair 
$(x,x+s-a)$ to the pair $(y,y+q-a)$. We have just proved 
that $\bar{g}$ leaves the set $\big\{ \, (x,x+s-a) : x \in \Z_n, s \in S \, \big\}$ fixed. As the latter set is 
the set of all directed edges of the digraph $\cay(\Z_n,S-a),$ $\bar{G} \le \aut(\cay(\Z_n,S-a))$.
For an automorphism $h$ of $\cay(\Z_n,S-a),$ define the permutation $g$ of $V$ by  
\[ g : x^\varepsilon \mapsto \begin{cases} 
(x^h)^+ & \text{ if } \varepsilon = +, \\
((x-a)^h+a)^- & \text{ if } \varepsilon = -,  
\end{cases} \quad x \in \Z_n, \; \varepsilon \in \{+,-\}. \]
The reader is invited to check that the above permutation $g$ is an automorphism of $\Gamma$. 
It is clear that $g \in G$ and $\bar{g} = h;$ we conclude that $\bar{G} = \aut(\cay(\Z_n,S-a))$. 

Now, the proposition follows along the following equivalences:
\begin{eqnarray*}
(1) & \iff & \text{Every bicyclic group of $\Gamma$ is conjugate to $C$ in $G$ } \\
    & \iff & \text{Every regular cyclic subgroup of $\bar{G}$ is conjugate to $\bar{C}$ in $\bar{G}$ }  \\ 
    & \iff & (2). 
\end{eqnarray*}   
The last equivalence is Theorem \ref{B}. \QED

\noindent 
{\bf Remark.} 
We want to remark that the equality $\aut(\Gamma)_{0^+} = \aut(\Gamma)_{a^-}$ does not hold in general.
For example, take $\Gamma$ as the incidence graph of the projective space $\mathrm{PG}(d,q)$ where $d \ge 2$ and 
$q$ is a prime power (i.e., $\Gamma$ is the bipartite graph with color classes identified by the points and the hyperplanes, 
respectively, and the edges are defined by the incidence relation of the space). It is well-known that the   
space $\mathrm{PG}(d,q)$ admits a cyclic group of automorphisms (called a Singer subgroup)
acting regularly on the points and the hyperplanes, respectively. This shows that $\Gamma$ is isomorphic to a cyclic 
Haar graph, and we may identify the set of points with $\Z_n^+,$ and the set of hyperplanes with $\Z_n^-,$ where 
$n = (q^d-1)/(q-1)$. The automorphism group $\aut(\Gamma) = \mathrm{P \Gamma L}(d+1,q) \rtimes \Z_2$. As 
$\mathrm{P\Gamma L}(d+1,q)$ acts inequivalently on the points and the hyperplanes, $\aut(\Gamma)_{0^+}$ cannot be 
equal to $\aut(\Gamma)_{a^-}$ for any $a \in \Z_n$.

\section{ Haar graphs $H(\Z_{2m},\{0,u,v,v+m\})$}

In this section we prove Theorem  \ref{MAIN1} for Haar graphs $H(\Z_{n},S)$ satisfying certain additional conditions.

\begin{thm}\label{MAIN2}
Let $n = 2m$ and $S =\{0,u,v,v+m\}$ such that 
\begin{enumerate}[(a)]
\item $\Z_n = \langle u,v \rangle;$
\item $1< u < m,$ $ u \mid m;$ 
\item $\aut(H(\Z_n,S))_{0^+}$ leaves the set $\{0^-,u^-\}$ setwise fixed.
\end{enumerate}
Then $H(\Z_n,S) \cong H(\Z_n,T)$ if and only if there exist $a \in \Z_n^*$ and $b \in \Z_n$ such that 
\begin{enumerate}[(1)]
\item $a T+ b = S;$ or 
\item $a T + b = \{0,u+m,v,v+m\},$ and 
$2 \mid u,$ $2u \mid m,$ $u/2 \not\equiv v+ m/(2u) (\mathrm{mod\;}m/u)$.  
\end{enumerate}
\end{thm}

By (b) of Theorem \ref{MAIN2} we have $2u \le m$. We prove the extremal case, when $2u  = m$, separately. 
Notice that, in this case the conditions in (2) of Theorem \ref{MAIN2} that $2 \mid u,$ $2u \mid m$ and $u/2 \not\equiv v+ m/(2u) (\mathrm{mod\;}m/u)$
can be replaced by one condition: $u \equiv 2(\mathrm{mod\;}4)$.

\begin{lem}\label{Lem2}
Let $S$ be the set defined in Theorem \ref{MAIN2}. If  $2u = m,$ then $H(\Z_n,S) \cong H(\Z_n,T)$ if and only if 
there exist $a \in \Z_n^*$ and $b \in \Z_n$ such that
\begin{enumerate}[(1)]
\item $a T + b = S;$ or  
\item $a T + b = \{0,u+m,v,v+m\}$ and $u \equiv 2(\mathrm{mod\;}4)$.  
\end{enumerate}
\end{lem}

\proof 
Let $d = \gcd(n,v)$. Because of  $\sg{u,v} = \Z_n$  we have that $\gcd(u,v,n) = 1,$ i.e., $\gcd(n/4,v) = 1,$ and 
this  gives that $d \in \{1,2,4\}$. Note that, if $d \ne 1,$ then necessarily $2 \nmid \,  u$.
Then we can write $v = v_1 d,$ where $\gcd(v_1,n)=1$. Let $v_1^{-1}$ denote the inverse of $v_1$ in the group 
$\Z_n^*$. Then the following hold in $\Z_n$ (here use that $u = n/4$): 
\[
v_1^{-1} v = d,  \;  v_1^{-1} (v+2u) = d + 2u \text{ \ and \ } v_1^{-1} u \in \{u,3u\}. 
\] 
We conclude that $S$ is affinely equivalent to one of the sets $S_i(d),$ $i \in \{1,2\}$ and $d \in \{1,2,4\},$ where  
\[  S_1(d) = \{0,u,d,d+2u\} \text{ or }  S_2(d) = \{0,3u,d,d+2u\}. \]

The lemma follows from the following claims:
\begin{enumerate}[(i)]
\item $H(\Z_n,S_1(1)) \cong H(\Z_n,S_2(1))$.
\item $H(\Z_n,S_1(1)) \conga H(\Z_n,S_1(d))$ for $d \in\{2,4\};$
\item $H(\Z_n,S_1(d)) \conga H(\Z_n,S_2(d)) \iff d \in \{2,4\} \text{ or } (d=1 \text{ and }  u \not\equiv 2(\mathrm{mod\;}4));$ 
\end{enumerate}  

(i): Define the mapping $f : V \mapsto V$ by 
\[  f : x^\varepsilon \mapsto \begin{cases} 
x^\varepsilon & \text{ if } x \in \{0,1,\ldots.u-1\} \cup \{2u,\ldots.3u-1\},  \\ 
(x+2u)^\varepsilon & \text{ otherwise}. \end{cases}\]
We leave for the reader to verify that $f$ is in fact an isomorphism from $H(\Z_n,S_1(1))$ to $H(\Z_n,S_2(1))$ (compare the  
graphs in Figure 1;  here  the white vertices represent the color class $\Z_n^+,$ while the black ones represent the color class $\Z_n^-$).

\begin{figure}[t!] 
\centering

\begin{tikzpicture}[scale=1.25] ---
\foreach \i in {0,2,4,6} 
{
\draw (\i,0) circle (1.5pt) ;
\draw (\i,0.5) circle (1.5pt);
\draw (\i,1) circle (1.5pt);
\draw (\i,1.5) circle (1.5pt); 
}
\foreach \j in {0.04,4.04,6.04}
{
\draw (\j,0) -- (\j+1,0) (\j,0) -- (\j+1,0.5);
\draw (\j,0.5) -- (\j+1,0.5) (\j,0.5) -- (\j+1,0);
\draw (\j,1) -- (\j+1,1) (\j,1) -- (\j+1,1.5);
\draw (\j,1.5) -- (\j+1,1.5) (\j,1.5) -- (\j+1,1);
}
\draw (2.04,0) -- (3,1) (2.04,0) -- (3,1.5);
\draw (2.04,0.5) -- (3,1) (2.04,0.5) -- (3,1.5);
\draw (2.04,1) -- (3,0) (2.04,1) -- (3,0.5);
\draw (2.04,1.5) -- (3,0) (2.04,1.5) -- (3,0.5);

\foreach \k in {1,3,5}
{
\fill (\k,0) circle (1.5pt) (7,0) circle (1.5pt);
\fill (\k,0.5) circle (1.5pt) (7,0.5) circle (1.5pt);
\fill (\k,1) circle (1.5pt) (7,1) circle (1.5pt);
\fill (\k,1.5) circle (1.5pt) (7,1.5) circle (1.5pt);
\draw (\k,0) -- (\k+0.96,0) (\k,0) -- (\k+0.96,1);
\draw (\k,0.5) -- (\k+0.96,0.5) (\k,0.5) -- (\k+0.96,1.5);
\draw (\k,1) -- (\k+0.96,1) (\k,1) -- (\k+0.96,0.5);
\draw (\k,1.5) -- (\k+0.96,1.5) (\k,1.5) -- (\k+0.96,0);
}

\draw (0,0) node[above] {\footnotesize n-2} (0,0.5) node[above] {\footnotesize 2u-2};
\draw (1,0) node[above] {\footnotesize n-1} (1,0.5) node[above] {\footnotesize 2u-1};
\draw (2,0) node[above] {\footnotesize n-1} (2,0.5) node[above] {\footnotesize 2u-1};
\draw (3,0) node[above] {\footnotesize 3u} (3,0.5) node[above] {\footnotesize u};
\draw (4,0) node[above] {\footnotesize 3u} (4,0.5) node[above] {\footnotesize u} ;
\draw (5,0) node[above] {\footnotesize 3u+1} (5,0.5) node[above] {\footnotesize u+1} ;
\draw (6,0) node[above] {\footnotesize 3u+1} (6,0.5) node[above] {\footnotesize u+1} ;
\draw (7,0) node[above] {\footnotesize 3u+2} (7,0.5) node[above] {\footnotesize u+2} ;

\draw (0,1) node[above] {\footnotesize 3u-2} (0,1.5) node[above] {\footnotesize u-2};
\draw (1,1) node[above] {\footnotesize 3u-1} (1,1.5) node[above] {\footnotesize u-1};
\draw (2,1) node[above] {\footnotesize 3u-1} (2,1.5) node[above] {\footnotesize u-1};
\draw (3,1) node[above] {\footnotesize 2u} (3,1.5) node[above] {\footnotesize 0};
\draw (4,1) node[above] {\footnotesize 2u} (4,1.5) node[above] {\footnotesize 0} ;
\draw (5,1) node[above] {\footnotesize 2u+1} (5,1.5) node[above] {\footnotesize 1} ;
\draw (6,1) node[above] {\footnotesize 2u+1} (6,1.5) node[above] {\footnotesize 1} ;
\draw (7,1) node[above] {\footnotesize 2u+2} (7,1.5) node[above] {\footnotesize 2} ;
\end{tikzpicture}

\bigskip 

\begin{tikzpicture}[scale=1.25] ---
\foreach \i in {0,2,4,6} 
{
\draw (\i,0) circle (1.5pt) ;
\draw (\i,0.5) circle (1.5pt);
\draw (\i,1) circle (1.5pt);
\draw (\i,1.5) circle (1.5pt); 
}
\foreach \j in {0.04,4.04,6.04}
{
\draw (\j,0) -- (\j+1,0) (\j,0) -- (\j+1,0.5);
\draw (\j,0.5) -- (\j+1,0.5) (\j,0.5) -- (\j+1,0);
\draw (\j,1) -- (\j+1,1) (\j,1) -- (\j+1,1.5);
\draw (\j,1.5) -- (\j+1,1.5) (\j,1.5) -- (\j+1,1);
}
\draw (2.04,0) -- (3,1) (2.04,0) -- (3,1.5);
\draw (2.04,0.5) -- (3,1) (2.04,0.5) -- (3,1.5);
\draw (2.04,1) -- (3,0) (2.04,1) -- (3,0.5);
\draw (2.04,1.5) -- (3,0) (2.04,1.5) -- (3,0.5);

\foreach \k in {1,3,5}
{
\fill (\k,0) circle (1.5pt) (7,0) circle (1.5pt);
\fill (\k,0.5) circle (1.5pt) (7,0.5) circle (1.5pt);
\fill (\k,1) circle (1.5pt) (7,1) circle (1.5pt);
\fill (\k,1.5) circle (1.5pt) (7,1.5) circle (1.5pt);
\draw (\k,0) -- (\k+0.96,0) (\k,0) -- (\k+0.96,1.5);
\draw (\k,0.5) -- (\k+0.96,0.5) (\k,0.5) -- (\k+0.96,1);
\draw (\k,1) -- (\k+0.96,1) (\k,1) -- (\k+0.96,0);
\draw (\k,1.5) -- (\k+0.96,1.5) (\k,1.5) -- (\k+0.96,0.5);
}

\draw (0,0) node[above] {\footnotesize n-2} (0,0.5) node[above] {\footnotesize 2u-2};
\draw (1,0) node[above] {\footnotesize n-1} (1,0.5) node[above] {\footnotesize 2u-1};
\draw (2,0) node[above] {\footnotesize n-1} (2,0.5) node[above] {\footnotesize 2u-1};
\draw (3,0) node[above] {\footnotesize 3u} (3,0.5) node[above] {\footnotesize u};
\draw (4,0) node[above] {\footnotesize 3u} (4,0.5) node[above] {\footnotesize u} ;
\draw (5,0) node[above] {\footnotesize 3u+1} (5,0.5) node[above] {\footnotesize u+1} ;
\draw (6,0) node[above] {\footnotesize 3u+1} (6,0.5) node[above] {\footnotesize u+1} ;
\draw (7,0) node[above] {\footnotesize 3u+2} (7,0.5) node[above] {\footnotesize u+2} ;

\draw (0,1) node[above] {\footnotesize 3u-2} (0,1.5) node[above] {\footnotesize u-2};
\draw (1,1) node[above] {\footnotesize 3u-1} (1,1.5) node[above] {\footnotesize u-1};
\draw (2,1) node[above] {\footnotesize 3u-1} (2,1.5) node[above] {\footnotesize u-1};
\draw (3,1) node[above] {\footnotesize 2u} (3,1.5) node[above] {\footnotesize 0};
\draw (4,1) node[above] {\footnotesize 2u} (4,1.5) node[above] {\footnotesize 0} ;
\draw (5,1) node[above] {\footnotesize 2u+1} (5,1.5) node[above] {\footnotesize 1} ;
\draw (6,1) node[above] {\footnotesize 2u+1} (6,1.5) node[above] {\footnotesize 1} ;
\draw (7,1) node[above] {\footnotesize 2u+2} (7,1.5) node[above] {\footnotesize 2} ;
\end{tikzpicture}

\caption{ Haar graphs $H(\Z_n,S_1(1))$ and $H(\Z_n,S_2(1))$.}
\end{figure}
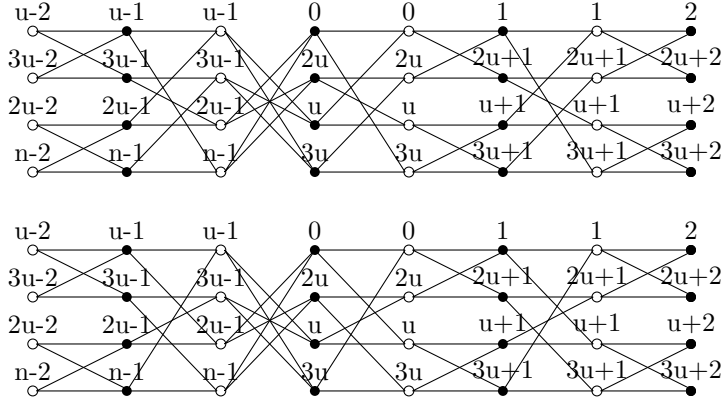

\medskip

(ii): Since $d \in \{2,4\},$ $u$ is an odd number. For $d \in \{2,4\}$ define $r_d \in \Z_n^*$ as follows:
\[ r_2 = \begin{cases} 
2 + u & \text{ if } u \equiv 1(\mathrm{mod\;}4), \\ 
2 + 3u & \text{ if } u \equiv 3(\mathrm{mod\;}4), \\ 
\end{cases} \; 
r_4 = \begin{cases} 
4 + u & \text{ if } u   \equiv 3(\mathrm{mod\;}4), \\ 
4 + 3u & \text{ if } u \equiv 1(\mathrm{mod\;}4). \\ 
\end{cases}
\] 
It can be directly checked that $r_d S_1(1) +  u = S_1(d),$ so 
$H(\Z_n,S_1(1)) \conga H(\Z_n,S_1(d))$ for $d \in\{2,4\}$.

\medskip

(iii): If $u$ is odd, then $(2u+1) S_1(d) = S_2(d),$ hence $H(\Z_n,S_1(d)) \conga H(\Z_n,S_2(d))$. 
Since $u$ is odd whenever $d \in \{2,4\},$ we are left with the case that $d = 1$ and $u$ is even.
If also $u \equiv 0(\text{mod } 4),$ then $(u+1) S_1(1) + 3u = S_2(1),$ and again $H(\Z_n,S_1(1)) \conga H(\Z_n,S_2(1))$.

Suppose that $d=1$ and $u \equiv 2(\mathrm{mod\;}4)$.  
We finish the proof by showing that in this case $H(\Z_n,S_1(1)) \not\conga H(\Z_n,S_2(1))$.
Suppose that, there is an affine transformation $\psi : x \mapsto r x + s,$ $r \in \Z_n^*$ and $s \in \Z_n,$ 
which maps the set $S_1(1)$ to $S_1(2)$. Then $1^\psi - (1+2u)^\psi = 2u$ in $\Z_n$. This implies that  
$\{1,1+2u\}^\psi = \{1,1+2u\}$ and $\{0,u\}^\psi = \{0,3u\},$ and hence     
\[     r + s \in \{1,1+2u\}  \text{ and }  r \{0,u\} + s = \{0,3u\}. \]
A direct analysis shows  that the above equations cannot hold if $u \equiv 2(\mathrm{mod\;}4)$. 
Thus $H(\Z_n,S_1(1)) \not\conga H(\Z_n,S_2(1)),$ this completes the proof of (iii). 
\QED 

\medskip

Now, we turn to the case when $2u \ne m$.  Recall that the canonical bicyclic group $C$ is generated by the 
permutation $c$ defined in (1). For a divisor $\ell \mid n,$ $C^\ell$ will denote the subgroup of  $C$ generated by $c^\ell$.
It will be convenient to denote by $\delta_{\ell}$ the partition of $V$ into the orbits of $C^\ell,$ 
i.e., $\delta_\ell = \orb(C^\ell,V)$. Furthermore, we set $\eta_{n,\ell}$ for the homomorphism $\eta_{n,\ell} : \Z_n \to \Z_\ell$ defined by $\eta_{n,\ell}(1)=1$. 

Observe that, if $\delta_\ell$ is, in addition, a  system of blocks for the group $A = \aut(H(\Z_n,S)),$ then   
an action of $A$ can be defined on $H(\Z_\ell,\eta_{n,\ell}(S))$ by letting $g \in A$ act as 
\begin{equation}\label{act}
(x^\varepsilon)^g = y^\varepsilon  \iff 
\big\{ z^-  : z \in \eta_{n,\ell}^{-1}(x)\big\}^g = \big\{ z^-  : z \in \eta_{n,\ell}^{-1}(y)\big\}, \;
x \in \Z_\ell, \; \varepsilon \in \{+,-\}. 
\end{equation}
We denote the corresponding kernel by $A_{(\delta_\ell)},$ and the image of $g \in G$ by $g^{\delta_\ell}$. 
Note that, if $X$ is a bicyclic group of $H(\Z_n,S),$ then 
$X^{\delta_\ell} = \{x^{\delta_\ell} : x \in X\}$ is a bicyclic group of $H(\Z_\ell,\eta_{n,\ell}(S))$.

\medskip

Let $S=\{0,u,v,v+m\}$ be the subset of $\Z_n$ defined in Theorem \ref{MAIN2}. Let $\delta$ be the partition of $V$ defined by  
\begin{equation}\label{D}
\delta = \big\{ X\cup X^{\psi_{1,0,0}} : X \in \orb(C^u,V) \big\},
\end{equation}
where $\psi_{1,0,0}$ is defined in \eqref{varphi psi}. 
We write $\delta = \{ V_0, \ldots, V_{u-1} \},$ where 
\[V_i = \big\{  (i v + j u)^+,(i v+ j u)^-   : j \in \{0,1,\ldots,(n/u)-1\} \big\}.\]

\begin{figure}
\centering

\begin{tikzpicture}[scale=1.5] ---


\foreach \i in {0,3,6} 
{
\fill (\i,3) circle (1.5pt) ;
\fill (\i,2.5) circle (1.5pt);
\fill (\i,2) circle (1.5pt) ;
\fill (\i,1.5) circle (1.5pt);
\fill (\i,1) circle (1.5pt) ;
\fill (\i,0.5) circle (1.5pt);
\draw[line width=.75 pt] (\i,3) -- (\i+.1,2.8);
\draw[line width=.75 pt] (\i,2.5) -- (\i+.3,2.3);
\draw (\i,.2) node {\footnotesize $\vdots$};
}


\foreach \j in {1.5,4.5,7.5} 
{
\draw (\j,3) circle (1.5pt); \draw[line width=1pt] (\j-0.04,3) -- (\j-1.5,3)  (\j-0.04,3) -- (\j-1.5,2);
\draw (\j,2.5) circle (1.5pt); \draw[line width=1pt] (\j-0.04,2.5) -- (\j-1.5,2.5)  (\j-0.04,2.5) -- (\j-1.5,1.5);
\draw (\j,2) circle (1.5pt); \draw[line width=1pt] (\j-0.04,2) -- (\j-1.5,2)  (\j-0.04,2) -- (\j-1.5,1) ;
\draw (\j,1.5) circle (1.5pt); \draw[line width=1pt] (\j-0.04,1.5) -- (\j-1.5,1.5)  (\j-0.04,1.5) -- (\j-1.5,.5);
\draw (\j,1) circle (1.5pt); \draw[line width=1pt]  (\j-0.04,1) -- (\j-1.5,1) (\j-0.04,1) -- (\j-.8,0.4);
\draw (\j,0.5) circle (1.5pt); \draw[line width=1pt] (\j-0.04,.5) -- (\j-1.5,0.5) (\j-0.04,0.5) -- (\j-.4,.2);
\draw (\j,.2) node {\footnotesize $\vdots$};
}


\foreach \k in {1.5,4.5}
{
\draw (\k+0.04,0.5) -- (\k+1.5,0.5) (\k+0.04,0.5) -- (\k+1.5,1);
\draw (\k+0.04,1) -- (\k+1.5,1) (\k+0.04,1) -- (\k+1.5,.5);
\draw (\k+0.04,1.5) -- (\k+1.5,1.5) (\k+0.04,1.5) -- (\k+1.5,2);
\draw (\k+0.04,2) -- (\k+1.5,2) (\k+0.04,2) -- (\k+1.5,1.5);
\draw (\k+0.04,2.5) -- (\k+1.5,2.5) (\k+0.04,2.5) -- (\k+1.5,3);
\draw (\k+0.04,3) -- (\k+1.5,3) (\k+0.04,3) -- (\k+1.5,2.5);
}


\draw (0,3) node[above] {\small $0$};   \draw (1.5,3) node[above] {\small $0$}; 
\draw (3,3) node[above] {\small $v$};   \draw (4.5,3) node[above] {\small $v$};
\draw (6,3) node[above] {\small $2v$}; \draw (7.5,3) node[above] {\small $2v$};

\draw (0,2.5) node[above] {\small $m$};       \draw (1.5,2.5) node[above] {\small $m$};
\draw (3,2.5) node[above] {\small $v+m$};   \draw (4.5,2.5) node[above] {\small $v+m$};
\draw (6,2.5) node[above] {\small $2v+m$}; \draw (7.5,2.5) node[above] {\small $2v+m$};

\draw (0,2) node[above] {\small $u$};       \draw (1.5,2) node[above] {\small $u$};
\draw (3,2) node[above] {\small $v+u$};   \draw (4.5,2) node[above] {\small $v+u$};
\draw (6,2) node[above] {\small $2v+u$}; \draw (7.5,2) node[above] {\small $2v+u$};

\draw (0,1.5) node[above] {\small $u+m$};       \draw (1.5,1.5) node[above] {\small $u+m$};
\draw (3,1.5) node[above] {\small $v+u+m$};   \draw (4.5,1.5) node[above] {\small $v+u+m$};
\draw (6,1.5) node[above] {\small $2v+u+m$}; \draw (7.5,1.5) node[above] {\small $2v+u+m$};

\draw (0,1) node[above] {\small $2u$};       \draw (1.5,1) node[above] {\small $2u$};
\draw (3,1) node[above] {\small $v+2u$};   \draw (4.5,1) node[above] {\small $v+2u$};
\draw (6,1) node[above] {\small $2v+2u$}; \draw (7.5,1) node[above] {\small $2v+2u$};

\draw (0,0.5) node[above] {\small $2u+m$};       \draw (1.5,0.5) node[above] {\small $2u+m$};
\draw (3,0.5) node[above] {\small $v+2u+m$};   \draw (4.5,0.5) node[above] {\small $v+2u+m$};
\draw (6,0.5) node[above] {\small $2v+2u+m$}; \draw (7.5,0.5) node[above] {\small $2v+2u+m$};

\draw (.75,-.25)   node {$V_0$};
\draw (3.75,-.25) node {$V_1$};
\draw (6.75,-.25) node {$V_2$};
\end{tikzpicture}

\caption{ The Haar graph $H(\Z_n,S)$ } 
\end{figure}
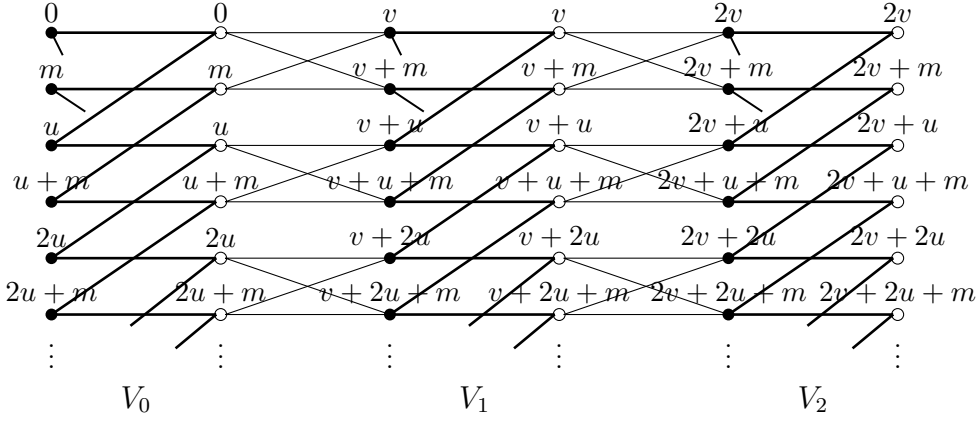

A part of $H(\Z_n,S)$  is drawn  in  Figure 2 using the partition $\delta$. White and black colors represent 
again the color classes $\Z_n^+$ and $\Z_n^-,$ respectively.
For $i \in \{0,1,\ldots,u-1\},$ let $e_i$ be the involution of $V$  defined by
\[ e_i : x^\varepsilon \mapsto \begin{cases} 
(x+m)^\varepsilon & \text{ if } x^\varepsilon \in V_i, \\ 
x^\varepsilon & \text{otherwise}.\end{cases} \]
It is clear that each $e_i \in \aut(H(\Z_n,S)),$ and also that $e_i e_j = e_j e_i$ for all $i,j \in \{0,1,\ldots,u-1\}$.  
Let $E = \sg{e_0,e_1,\ldots,e_{u-1}}$. Thus $E \le \aut(H(\Z_n,S))$ and $E \cong \Z_2^u$.
For a subset $I \subseteq \{0,1,\dots,u-1\}$ let $e_I$ be the element in $E$ defined by  
$e_I = \prod_{i \in I}e_i$.

\medskip

The following lemma about imprimitivity systems of blocks will be used throughout the paper.  

\begin{lem}\label{Lemnew}
Let $\Gamma = H(\Z_n,R)$ be a Haar graph and suppose that $T \subset R$ such that the point stabilizer $\aut(\Gamma)_{0^+}$ fixes setwise $T^-$, and let $d = |\langle T - T \rangle|,$ where 
$T - T = \{ t_1 - t_2 : t_1,t_2 \in T\}$.  The the partition $\pi$ of $V$ defined by 
$$ \pi = \{  X \cup X^{\psi_{1,t,-t}} : X \in \orb(C^{n/d},V) \}, \text{ where } t \in T $$
is a system of blocks for $\aut(\Gamma)$.\footnote{Notice that, $\pi$ does not depend of the choice of the element $t \in T$.} 
\end{lem}

\proof  Put $A = \aut(\Gamma)$.
Since $T$ is fixed setwise by $A_{0^+}$, we may write 
$$ T = T_1 \cup \cdots \cup T_k,$$ 
where $T_i^-$ is an $A_{0^+}$-orbit for every $i \in \{1,2,\dots,k\}$.
Choose an arc $(0^+,t_i^-)$ of $\Gamma$ where we fix an element $t_i \in T_i$ for every 
$i \in \{1,\dots,k\}$. 
We claim that, the orbital graph of $A$ containing $(0^+,t_i^-)$ is self-paired, and in fact it is 
equal to the Haar graph $H(\Z_n,T_i)$ (for a definition of an orbital graph, see \cite{DixM96}). 
 
Define $\bar{A}$ as the color preserving subgroup  of $A$. 
Then $A = \bar{A} \rtimes \langle \psi_{-1,0,0} \rangle$. Also, $\bar{A} = A_{0^+} C,$ as $C$ is transitive on $\Z_n^+$. 
Then the orbit of the arc $(0^+,t_i^-)$ under $A$ is 
\begin{eqnarray*}
(0^+,t_i^-)^A & = & (0^+,t_i^-)^{A_{0^+} C \langle \psi_{-1,0,0} \rangle} = 
\{ (0^+,t_i'^{\, -}) : t_i' \in T_i \}^{C \langle \psi_{-1,0,0} \rangle} \\
& = & \{ (j^+,(j+t_i')^{\, -}) : t_i' \in T_i, j \in \Z_n  \}^{\langle \psi_{-1,0,0} \rangle} \\ 
& = & \{ (j^+,(j+t_i')^{\, -}) : t_i' \in T_i, j \in \Z_n  \} \cup \{ ((-j)^-,(-j-t_i')^{\, +}) : t_i' \in T_i, j \in \Z_n  \} , \\
& = & \{ (j^+,(j+t_i')^{\, -}) : t_i' \in T_i, j \in \Z_n  \} \cup \{ ((j+t_i')^-,j^+) : t_i' \in T_i, j \in \Z_n  \} , 
\end{eqnarray*}
which is clearly equal to the set of the set of arcs of $H(\Z_n,T_i)$. The claim is proved.

Since $H(\Z_n,T_i)$ is an orbital graph, $A \le \aut(H(\Z_n,T_i))$. Combining this with 
$H(\Z_n,T) = \cup_{i=1}^k H(\Z_n,T_i),$ we have that $A \le \aut(H(\Z_n,T)$. 
Let $\Sigma$ be the connected component of $H(\Z_n,T)$ which contains $0^+$. Obviously, the 
set $W$ of vertices contained in $\Sigma$ is a block for $A$. 
It is easy to verify that $W =  X \cup X^{\psi_{1,t,-t}}$ where $X$ is the orbit of $0^+$ under $C^{n/d}$. 
The lemma follows.  \QED

\begin{lem}\label{Lem3}
Let $S$ be the set defined in Theorem \ref{MAIN2}. If  $2u \ne m,$ 
then the stabilizer $\aut(H(\Z_n,S))_{0^+}$ is given as follows. 
\begin{enumerate}[(1)]
\item If $u \not\equiv 2v(\mathrm{mod\;}m/u),$ then $\aut(H(\Z_n,S))_{0^+} = E_{0^+}$.
\item If $u \equiv 2v(\mathrm{mod\;}m/u),$ then $\aut(H(\Z_n,S))_{0^+} = E_{0^+} \times F$ for a subgroup 
$F \le \aut(H(\Z_n,$ $S))_{0^+},$ $|F| =2$.
\end{enumerate}
\end{lem}

\proof For short we put $\Gamma = H(\Z_n,S)$ and $A = \aut(\Gamma)$. 
Consider the partition $\delta$ defined in \eqref{D}. Applying Lemma \ref{Lemnew} with $R=S$, $T=\{0,u\}$ 
and $t=0$, we obtain that $\delta$  is a system of blocks for $A$.
The quotient graph $\Gamma/\delta$ is a $u$-circuit if $u > 2$ and 
a $2$-path if $u=2$. Let $g \in A_{0^+}$. Then $g$ fixes the directed edge $(V_0,V_1)$ of $\Gamma/\delta,$ 
hence it must fix all sets $V_i$. Thus $A_{0^+} \le A_{(\delta)},$ where $A_{(\delta)}$ is the kernel of the action of $A$ on 
$\delta$.  

Consider the action of $A_{0^+}$ on $V_0$. The  corresponding kernel is 
$A_{(V_0)},$  the pointwise stabilizer of $V_0$ in $A$, and the corresponding image 
is a subgroup of $\aut(\Gamma[V_0]),$ where $\Gamma[V_0]$ is the subgraph of $\Gamma$ induced by $V_0$.  
Using that $2u \ne m$, we show next that $A_{(V_0)} = E_{0^+}$.
 It is clear that $A_{(V_0)} \geq E_{0^+}$. 
We are going to prove that $A_{(V_0)} \le E_{0^+}$ also holds. 
Let $g \in A_{(V_0)}$. Then for a suitable element $e \in \langle e_1 \rangle$ the product 
$g e$ fixes pointwise $V_0$ and fix the vertex $v^+$ from block $V_1$. 
Thus $g e$ acts on $V_1$  as the identity or the unique reflection of the circuit $\Gamma[V_1]$ that fixes 
$v^+$. If this action is not the identity, then $g e$ switches $v^-$ and $(v+n-u)^-,$ and so it must switch $(v+u)^+$ and 
$(v+n-u)^+$. This implies that $(v+n-u)^+ = (v+u+m)^+,$ contradicting that $2u \ne m$. We obtain that 
$g e$ acts as the identity also on $V_1$.  Continuing in tis way we find that $g e'$ is the identity with a suitable choice of 
$e' \in E_{0^+},$ hence $g = e'$.

The equality $A_{(V_0)} = E_{0^+}$ together with $\aut(\Gamma[V_0]) \cong D_{4u}$ imply that 
$|A_{0^+} : E_{0^+}| \le 2$. Moreover,  $|A_{0^+} : E_{0^+}| = 2$  holds exactly when $A_{0^+}$ contains an involution $g$ for which  $g : 0^- \leftrightarrow  u^-$. In the latter case  $A_{0^+} = E_{0^+} \times \sg{g}$ as $g$ centralizes $E$ (observe that $g$ is in the kernel $A_{( \delta )},$ hence acts on a block $V_i$ as an element of 
$D_{2 n/u},$ while $E$ acts on $V_i$ as the center $Z(D_{2n/u})$.) 
We  settle the lemma by proving the following equivalence :
\begin{equation}\label{iff}
A_{0^+} \cong E_{0^+} \times \Z_2  \iff u \equiv 2v(\text{mod }m/u).
\end{equation}
 
Suppose first that $A_{0^+} = E_{0^+} \times \sg{g},$ where $g \in A_{0^+}$ and $g : 0^- \leftrightarrow  u^-$.
By (c) of Theorem 3.1, $\{v^-,(v+m)^-\}^{A_{0^+}} = \{v^-,(v+m)^-\}$.
Applying Lemma \ref{Lemnew} with $R=S$, $T=\{v,v+m\}$ 
and $t=v$, we obtain that the set $B = \{0^+,m^+,v^-,(v+m)^-\}$ is a block for $A$. 
The induced graph $\Gamma[B]$ is a $4$-circuit. Denote by $A_{\{B\}}$ the setwise stabilizer of $B$ in $A,$ and 
by $A_{\{B\}}^B$ the permutation group of $B$ induced by $A_{\{B\}}$.  As $\Gamma[B]$ is a $4$-circuit, 
$A_{\{B\}}^B \le D_8$. This gives that $\{0^+,m^+\}$ is a block for $A_{\{B\}}^B,$ and therefore it is also a 
block for $A$. We conclude that  
$\delta_m = \{ X : X \in \orb(C^m,V)\}$ is a system of blocks for $A$. 
Consider the action of $A$ on $H(\Z_m,\eta_{n,m}(S))$ defined in \eqref{act}.
Then $E \le A_{(\delta_m)}$, while $g \notin A_{(\delta_m)}$. 
This implies that $g^{\delta_m}$ is an automorphism of 
$H(\Z_m,\eta_{n,m}(S))$ which normalizes its canonical bicyclic group. 
This means that $g^{\delta_m} = \varphi_{r,s,t}$ for some $r \in \Z_m^*$ and $s,t \in \Z_m$. 
Using that $g^{\delta_m} : 0^+ \mapsto 0^+$ and $0^- \mapsto \eta_{n,m}(u)^-,$ we find that $t=0$ and 
$s = \eta_{n,m}(u),$ and so 
\begin{equation}\label{Ad}
A^{\delta_m} = \sg{D^{\delta_m},\varphi_{r,\eta_{n,m}(u),0}}.
\end{equation}
Also,  $g^{\delta_m} :  \eta_{n,m}(u)^- \mapsto 0^- $ and $\eta_{n,m}(v)^- \mapsto \eta_{n,m}(v)^-,$ hence  
$ r \eta_{n,m}(u) =  - \eta_{n,m}(u)$  and $r \eta_{n,m}(v) = \eta_{n,m}(v-u)$ 
hold in $\Z_m$. From these $r \equiv -1(\mathrm{mod\;}m/u)$  and 
$rv \equiv v-u(\text{mod }m/u),$ i.e., $u \equiv 2v(\mathrm{mod\;}m/u)$.  
 The implication ``$\Rightarrow$''  in \eqref{iff}  is now proved.  

\medskip

Suppose next that $u \equiv 2v(\text{mod }m/u)$. Define the permutation $g$ of $V$ by    
\[ g : (vi+uj)^\varepsilon  \mapsto \begin{cases}  
\big( v i - (i+j)u \big)^+ &  \text{ if } \varepsilon = +, \\ 
\big( v i - (i+j-1)u \big)^- & \text{ if } \varepsilon = -,\end{cases}\]
where $ i \in \{0,1,\ldots,u-1\}$  and $j \in \{0,1,\ldots,n/u-1\}$. 
It can be verified directly that $g \in A_{0^+},$ ${0^+}^{\, g} = 0^+$ and 
$g : 0^- \leftrightarrow u^-,$ implying that $A_{0^+} = E_{0^+} \times \sg{g}$. Thus  part ``$\Leftarrow$''  of \eqref{iff} is 
also true. The lemma is proved.  
\QED 

\begin{lem}\label{Lem4}
Let $S$ be the set defined in Theorem \ref{MAIN2}. If  $2u \ne m,$ 
then for the normalizer $N_{\aut(H(\Z_n,S))}(C)$ of $C$ in $\aut(H(\Z_n,S)),$

\begin{equation}\label{N} 
\big| \aut(H(\Z_n,S)) : N_{\aut(H(\Z_n,S))}(C) \big| = \begin{cases} 
2^{u-2} & \text{ if } 2 \mid u \text{ and } \big( \, u \not\equiv 2v(\mathrm{mod\;}m/u) \text{ or } \\ 
              & \; u/2 \equiv v(\mathrm{mod\;}m/u) \, \big), \\  
2^{u-1} & \text{ otherwise}. \end{cases} 
\end{equation}
\end{lem}

\proof  
For short we put $A = \aut(H(\Z_n,S))$ and $N = N_A(C)$. 
Since $A = D A_{0^+}$ and $D \le N,$ $N = D (N \cap A_0^+) $. The cases (1)  and (2) in Lemma \ref{Lem3} 
 are considered separately.

\medskip 

\noindent {\sc Case 1.}  $u \not\equiv 2v(\mathrm{mod\;}m/u)$.  

\medskip 

In this case, from Lemma \ref{Lem3},  $A_{0^+} = E_{0^+},$ hence $|A| = 2^u n$. Let $g \in N \cap A_0^+$.  
Since $g \in E_{0^+},$ it follows quickly 
that $g=1_A$ or   $2 \mid u$ and $g = e_1e_3 \cdots  e_{u-1}$.  Combining this with 
$N = D (N \cap G_0^+) $ we find that  $|N| = 4n$ if $2 | u$, and 
$|N| = 2n$ if $2 \nmid u$.  Formula \eqref{N} follows.  
 
\medskip

\noindent {\sc Case 2.}  $u \equiv 2v(\text{mod }m/u)$.

 \medskip

From Lemma \ref{Lem3}, $A_{0^+}  =  E_{0^+} \times F$ for a subgroup $F \le A_{0^+},$ $|F|=2,$ hence 
$|A| = 2^{u+1} n$. 
It follows from the proof of Lemma \ref{Lem3} that, there exists $r \in \Z_m^*$ such that the following hold:  
\[ r \eta_{n,m}(u) = -\eta_{n,m}(u)  \text{ and } r \eta_{n,m}(v) = \eta_{n,m}(v-u). \]
Let $s \in \Z_n^*$ such that $\eta_{n,m}(s) = r$. Then  
\begin{equation}\label{su sv}
s u \in \{-u,-u+m\} \text{ and } s v \in \{v-u,v-u+m\}. 
\end{equation} 

Suppose that $2 \nmid  u$. Then we get as before that $N \cap E_{0^+}$ is trivial. 
Notice also that,  $u \equiv 2v + m/u(\mathrm{mod\;}n/u),$ which follows from the assumption that 
$u \equiv 2v(\mathrm{mod\;}m/u)$ and that $2 \nmid u$. Thus $2 \nmid m$ and $2 \mid (u+m),$ implying that 
in  \eqref{su sv} we have $s u = -u$. We obtain that $\varphi_{s,u,0} \in N \cap (A_{0^+} \setminus E_{0^+}),$ and so 
$|N \cap A_{0^+}| = 2$. 

Suppose next that $2 \mid u$. Then $|N \cap E_{0^+}| = 2$. It is easily 
seen that $|N \cap A_{0^+}| = 4$ if and only if there exists $r \in \Z_n^*$ such that 
$r u = -u $ and $r v = v-u$ hold in $\Z_n$. Consider the following system of linear congruences: 
\begin{equation}\label{system}
x u \equiv u(\mathrm{mod\;}n), \; x v \equiv v-u(\mathrm{mod\;}n).
\end{equation}
From the first congruence we can write $x$ in the form $x = y n/u-1$. Substitute this into the second congruence. We obtain that 
$y v n/u \equiv 2v-u(\mathrm{mod\;}n)$. This has a solution if and only if $\gcd(v n/u,n) \mid (2v-u)$. 
Suppose that $\gcd(v,n) \ne 1$. Using that $\sg{u,v} = \Z_n$ and that $ 2 \mid u,$ we obtain that $\gcd(v,m/u) \ne 1$. 
However, then from the assumption that 
$u \equiv 2v(\mathrm{mod\;}m/u)$ it follows that also $\gcd(v,u) \ne 1,$ which contradicts that $\sg{u,v} = \Z_n$.
Hence $\gcd(v,n)=1,$ $\gcd(v n/u,n) = n/u,$ and so  \eqref{system} has a solution if and only if 
$u \equiv 2v(\text{mod }n/u),$ or equivalently, 
$u/2 \equiv v(\mathrm{mod\;}m/u)$ (recall that $2 \mid u$ and $u \mid m$).  
It is not hard to show that any solution to \eqref{system} is necessarily  prime to $n,$ hence is in $\Z_n^*$. 
The above arguments can be summarized as follows:  $|N| = 8n$ if  
$2 \mid u$ and $u/2 \equiv v(\mathrm{mod\;}m/u),$ and $|N| = 4n$  otherwise.  This is consistent with \eqref{N}.  
The lemma is proved. 
\QED 
 
\begin{lem}\label{Lem5}
Let $r \in \Z_n^*,$ $r \ne 1$ and $s \in \Z_n$ such that the permutation  $\varphi_{r,s,0}$ is of order $2$. Then the  
group  $\sg{D,\varphi_{r,s,0} }$ contains a bicyclic subgroup different from $C$ if and only if 
$8 \mid n,$ $r = n/2+1,$ and $s = 0$ or $s = n/2$.
\end{lem}

\proof Suppose that $\sg{D,\varphi_{r,s,0} }$ contains a bicyclic subgroup $X$ such that $X \ne C$. 
Then $X$ is generated by a permutation in the form $c^i \varphi_{r,s,0}$. Since $\varphi_s^2 = id_V,$  
$r^2 = 1$ in $\Z_n,$ and we calculate that 
 $(c^i \varphi_{r,s,0})^2$ sends $x^+$ to $ (x+r(r+1)i)^{\, +}$ for every $x \in \Z_n$. 
That $\Z_n^+$ is an orbit of $X$ is equivalent to the condition that $\gcd(n,r+1) = 2$. Using this and that 
$r^2-1 = (r-1)(r+1) \equiv 0(\mathrm{mod\;}n),$ we find that $n/2$ divides $r-1$, so $r=1$ or $r=n/2 + 1$. Since $r \ne 1$, we have that $r = n/2 + 1$ and $8 \mid n$.
Then $(\varphi_{r,s,0})^2$ sends $x^-$ to $(x+(n/2+2)s)^{\, -}$. Since $(\varphi_{r,s,0})^2 = id_V,$ we obtain that 
$s = 0$ or $s = n/2$. 

On the other hand, it can be directly checked that,  if $8 \mid n,$ $r=n/2+1$ and $s \in \{0,n/2\},$ then the permutation 
$c \varphi_{r,s,0}$ generates a bicyclic subgroup of $\sg{D,\varphi_{r,s,0} }$. Obviously, this  bicyclic subgroup 
cannot be $C$. The lemma is proved. 
\QED 

Everything is prepared to prove the main result of the section.

\bigskip

\noindent{\sc Proof of Theorem \ref{MAIN2}.}  
The case that $2u = m$ is settled already  in Lemma \ref{Lem2}, hence let $2u \ne m$. 
We consider the action of $A = \aut(H(\Z_n,S))$ on the system of blocks $\delta_m$  defined in \eqref{act}. 
We claim that the corresponding image $A^{\delta_m}$ has a unique 
bicyclic subgroup (which is, of course, $C^{\delta_m}$). 

This is easy to see if $A_{0^+} = E_{0^+},$ because in this case $A^{\delta_m} = (D A _{0^+})^{\delta_m} = D^{\delta_m}$.

Let $A_{0^+} \ne E_{0^+}$. Then $A_{0^+} = E_{0^+} \times F$ for some subgroup $F,$ $|F|=2$. 
By \eqref{Ad}, $A^{\delta_m} = \sg{D^{\delta_m},\varphi_{r,\eta_{n,m}(u),0}}$.
Also, $r \equiv -1(\mathrm{mod\;}m/u),$ hence $r \ne 1$ in $\Z_m$. By Lemma 3.5,  
$A^{\delta_m}$ contains more than one bicyclic subgroup 
if and only if $8 \mid m,$ $r = m/2+1$ and  $\eta_{n,m}(u) \in \{0,m/2\}$. In the latter case $u \in \{m,m/2\},$ 
which is impossible as  $u < m/2$. Hence 
$A^{\delta_m}$ contains indeed a unique bicyclic subgroup.

We calculate below the number of bicyclic groups of $H(\Z_n,S);$ we denote this number by $\GBG$. 
Since $G = D A_{0^+},$  every $g \in G$  can be written as $g = x y$ with $x \in D$ and $y \in A_{0^+}$.
If $\sg{g}$ is a bicyclic group of $H(\Z_n,S),$ then $\sg{g}^{\delta_m}$ is a 
bicyclic subgroup of $A^{\delta_m},$ and so $\sg{g}^{\delta_m}  = C^{\delta_m}$.
This implies that  $x = c^i \in C$ for some $i \in \{1,\dots ,n-1\}$ with $\gcd(i,m)=1,$ and that 
$y \in E_{0^+}$.   
Then $y = e_I$ for a subset $I \subseteq \{1,\dots,u-1\}$. 
It follows  by induction that  
\[ (c^i e_I)^u = c^{ui} e_I e_{I+i} \cdots e_{I+(u-1)i}.\]
Here for $k \in \{1,\dots,u-1\},$ the set $I+ki = \{x+ki : x \in I\},$ where addition is taken modulo 
$u$. Since $\gcd(i,m)=1$ and $u \mid m,$ $\gcd(i,u)=1,$ from which 
\[ e_I e_{I+i} \cdots e_{I+(u-1)i} = (e_0e_1 \cdots e_{u-1})^{|I|} = c^{m|I|}.\]
Therefore, $(c^i e_I)^u = c^{\, u( \, i+ \frac{m}{u}|I| \,)}.$
This shows that $\sg{c^i e_I}$ is a bicyclic group 
if and only if $c^i e_I$  is of order $n,$  or equivalently, 
\begin{equation}\label{gcd}
\gcd\Big( i+\frac{m}{u} |I|,  \, \frac{2m}{u} \Big) = 1.
\end{equation}

Let $2 \mid (m/u)$. Then $2 \mid m$ and $i$ is odd. Hence 
\eqref{gcd} always holds. We obtain that the number of elements in $A$ which generate a bicyclic group is $\phi(n) 2^{u-1}$, where $\phi$ denotes the Euler's totient function. This implies that $\GBG = 2^{u-1}$.

Suppose next that $2 \nmid (m/u)$. Now, if $2 \mid u,$ then $2 \mid m,$ hence $2 \nmid i,$ and so \eqref{gcd} is true if and only if 
$|I|$ is even. We deduce from this that $\GBG = 2^{u-2}$. 

Finally, if $2 \nmid u,$ then $2 \nmid m,$ and in this case 
\eqref{gcd} holds if and only if $\gcd(i,n)=1$ and $|I|$ is even, or $\gcd(i,n)=2$ and 
$|I|$ is odd.  We calculate that $\GBG =  2^{u-1}$. All cases have been considered, our calculations are summarized in the formula:
\begin{equation}\label{GBG}
\GBG = \begin{cases} 
2^{u-2} & \mbox{ if } 2 \mid u \text{ and } 2 \nmid  (m/u),  \\ 
2^{u-1} & \mbox{ otherwise}.
\end{cases}
\end{equation}

Let $\Xi$ be a bicyclic base of $H(\Z_n,S)$. 
By \eqref{N} and \eqref{GBG} we obtain that,  $|\Xi| > 1$ if and only if  
\[ |A : N_A(C)| = 2^{u-2} \text{ and } \GBG = 2^{u-1}. \]
This happens exactly when 
\[  \big( \, 2 \mid u \text{ and } (u \not\equiv 2v(\mathrm{mod\;}m/u) \text{ or } u/2 \equiv v(\mathrm{mod\;}m/u)) \, 
\big) \text{ and } \big( 2 \nmid u \text{ or } 2u \mid m \big).\]
After some simplification, 
\[ |\Xi|  > 1 \iff 2 \mid u, \; 2u \mid m \text{ and }  u/2 \not\equiv v+ m/(2u) (\mathrm{mod\;}m/u). \]

Suppose that $|\Xi| > 1$. Then $A$ contains exactly $2^{n-1}$ bicyclic subgroups, $2^{n-2}$ of which are conjugate to 
$C$. These $2^{n-1}$ subgroups are enumerated as:  
$\sg{ c e _I}, \; I \subseteq \{1,\dots,u-1\}.$  For $i \in \{1,\dots,u-2\},$ $e_i c e_i = c e_{\{i,i+1\}}$. We can conclude that 
the bicyclic subgroups split into two conjugacy classes: 
\[ \big\{  \sg{c e_I} : I \subseteq \{1,\dots,u-1\}, |I| \text{ is even } \big\} \text{ and } \big\{  \sg{c e_I} : 
I \subseteq \{1,\dots,u-1\}, |I| \text{ is odd } \big\}. \] 
In particular, $|\Xi| = 2$. Choose $\xi$  from $\sym(V)$ which satisfies 
\[ \xi c \xi^{-1} = c e_1 \text{ and } \xi : 0^+ \mapsto 0^+, \; 0^- \mapsto 0^-. \] 
Then $\Xi$ can be chosen as $\Xi = \{id_V,\xi\}$. Also, $\{v^-,(v+m)^-\}^\xi = \{v^-,(v+m)^-\},$ and 
since $(ce_1)^{u+m}=c^u$, $(u^-)^\xi = (0^-)^{(ce_1)^{u+m} \xi} = (0^-)^{\xi c^{u+m}} = (u+m)^-$. 
 Thus $H(\Z_n,S)^{\, \xi} = H(\Z_n,\{0,u+m,v,v+m\})$. The theorem follows from Theorem \ref{Thm1}.
\QED 

\section{Proof of Theorem \ref{MAIN1}}

Theorem \ref{MAIN1} follows from Theorem \ref{MAIN2} and the following theorem.

\begin{thm}\label{MAIN3}
Let $H(\Z_n,S)$ be a connected  Haar graph such that $|S| = 4$ and $S$ is not a BCI-subset.
Then $n = 2m,$ and there exist $a \in \Z_n^*$ and $ b\in \Z_n$ such that   
$a S + b = \{0,u,v,v+m\}$ and  the conditions (a)-(c) in Theorem \ref{MAIN2} hold. 
\end{thm}

Before we prove Theorem\ref{MAIN3}  it is necessary to give three preparatory lemmas. 
For an element $i \in\Z_n,$ we denote by $|i|$ the order of $i$ viewed as an element   of the additive 
group $\Z_n$. Thus we have $|i| = n/\gcd(n,i)$.  

\begin{lem}\label{Lem6}
If  $R = \{i,n-i,j\}$ is  a generating subset of $\Z_n$ with $|i|$ odd, then $R$ is a CI-subset. 
\end{lem}

\proof For short we put $A = \aut(\cay(\Z_n,R))$ and denote by $A_0$ the stabilizer of $0 \in \Z_n$  in $A$. 
Clearly, $A_0$ leaves $R$ fixed. 
If $A_0$ acts on $R$ trivially, then $A \cong \Z_n,$ and the lemma follows by Theorem \ref{B}. 
If $A_0$ acts on $R$ transitively, then $\cay(\Z_n,R)$ is edge-transitive. This condition forces that 
$R$ is  a CI-subset (see \cite[page 320]{Li02}).

We are left with the case that $R$ consists of two orbits under $A_0$. 
These orbits must be $\{i,n-i\}$ and $\{j\}$. It is clear that $A_0$ leaves the subgroups 
$\sg{i}$ and $\sg{j}$ fixed; moreover, the latter set is fixed pointwise, and  
since  $|i|$  is odd, $\sg{i}$ consists of $(|i|-1)/2$ orbits under $A_0$, each of length $2,$ and one orbit of length $1$. 
We conclude that $\Z_n = \sg{i} \times \sg{j},$ and also that $A$ is permutation isomorphic to the permutation 
direct product 
$\big( (\Z_{|i|})_{\text right} \rtimes \sg{\pi} \big) \times  (\Z_{|j|})_{\text right}$. 
For $\ell \in \{|i|,|j|\},$ $(\Z_\ell)_{\text right}$ is generated by the affine transformation $x \mapsto x+1,$ 
and $\pi$ is the affine transformation $x \mapsto -x$.  
We leave for the reader to verify that the above group has a unique regular cyclic subgroup. 
The lemma follows by Theorem \ref{B}. \QED

\begin{lem}\label{Lem7}
Let $n=2m$ and $R = \{i,n-i,j,j+m\}$ be a subset of $\Z_n$ such that
\begin{enumerate}[(a)]
\item $R$ generates $\Z_n$;
\item $|i|$ is odd;
\item the stabilizer $\aut(\cay(\Z_n,R))_0$ leaves the set $\{i,n-i\}$ setwise fixed. 
\end{enumerate}
Then $R$ is a CI-subset. 
\end{lem}

\proof For short we put $A = \aut(\cay(\Z_n,R))$. Let $T$ be a subset of $\Z_n$ such that 
$\cay(\Z_n,R) \cong \cay(\Z_n,T)$ and 
let $f$ be an isomorphism from $\cay(\Z_n,R)$  to $\cay(\Z_n,T)$ such that $f(0) = 0$. 

Let us consider the subgraphs $\Gamma_1:= \cay(\Z_n,\{i,n-i\})$ and $\Gamma_2:=\cay(\Z_n,\{j,j+m\})$.
By condition $(c),$ the group $A$ preserves both of these subgraphs, that is, $A \le \aut(\Gamma_i)$ for $i \in \{1,2\}$.
Therefore, $f^{-1} (\Z_n)_{\text right} f \le A \le \aut(\Gamma_i),$ implying that $f$ maps $\Gamma_i$ to a  
Cayley graph $\cay(\Z_n,T_i)$ for both $i \in \{1,2\}$. Clearly, $T = T_1 \cup T_2$.  
Since $\Gamma_1= \cay(\Z_n,\{i,n-i\}) \cong \cay(\Z_n,T_1)$ it follows by a result of Sun \cite{Sun84} (see also \cite{Li02}), 
that $T_1 = a \{i,n-i\}$ for some $a \in \Z_n^*$. Letting $t_1 = a i,$ we have $T_1 = \{ t_1,n-t_1 \}$ such that 
$|i| = |t_1|$.  In the same way $T_2 = a' \{j,j+m\}$ for some $a' \in \Z_n^*,$ and letting $t_2 = a' j,$ we have  
$T_2 = \{t_2,t_2+m\}$ such that $|t_2| = |j|$.  
Since $f(0) = 0,$ $f$ maps $\{i,n-i\}$ to $T_1 = \{t_1,n-t_1\}$ and $\{j,j+m\}$ to $T_2 = \{t_2,t_2+m\}$. 

We claim that the partition of $\Z_n$ into the cosets of $\sg{m}$ is a system of blocks for $\aut(\Gamma_1),$ hence 
also for the group $A \le \aut(\Gamma_1)$.  Let us put $\bar{A} =  \aut(\Gamma_1)$. It easily seen that 
$\bar{A}_0$ is transitive on $\{j,j+m\}$. Therefore, by considering the action of $\bar{A}_0$ on the set $\{j,j+m\}$ we 
find that $\bar{A}_0 = \bar{A}_{0,j} \langle m_{\text right} \rangle = \langle m_{\text right} \rangle \bar{A}_{0,j},$ where 
$\bar{A}_{0,j} = \bar{A}_0 \cap \bar{A}_j,$ and $m_{\text right}$ is the permutation $x \mapsto x+m$ for every 
$x \in \Z_n$. From this $\bar{A}_0  \langle m_{\text right} \rangle = 
\langle m_{\text right} \rangle \bar{A}_0$. 
This means $\bar{A}_0  \langle m_{\text right} \rangle$ is a subgroup of $\bar{A}$ which clearly contains $\bar{A}_0$. 
It is well-known that the orbit $\orb_{\bar{A}_0  \langle m_{\text right} \rangle}(0)$ is a block for $\bar{A}$. Now the 
required statement follows as $\orb_{\bar{A}_0  \langle m_{\text right} \rangle}(0)  
= 0^{\bar{A}_0  \langle m_{\text right} \rangle} = 0^{\langle m_{\text right} \rangle} = \langle m \rangle$. 

Since the partition of $\Z_n$ into the cosets of $\sg{m}$ is a system of blocks for $A,$
the isomorphism $f$ induces an isomorphism from $\cay(\Z_m,\eta_{n,m}(R))$ to $\cay(\Z_m,\eta_{n,m}(T))$, we denote 
this isomorphism by $\bar{f}$. Note that, $\bar{f}(0)=0$ for the identity element $0 \in \Z_m$. 

The set $\eta_{n,m}(R)$ satisfies the conditions (a)-(c) of Lemma \ref{Lem6}, hence it is a CI-subset. 
This means that  $\bar{f}$ is equal to a permutation $x \mapsto rx$ for some $r\in \Z_m^*$. 
Let $s \in \Z_n^*$ such that $\eta_{n,m}(s)=r$. Then 
$\eta_{n,m}(si) = \eta_{n,m}(s)\eta_{n,m}(i)=\eta_{n,m}(t_1),$ and so the following holds in $\Z_n$: 
\begin{equation}\label{CI}
s i = t_1 \text{ or } s i =  t_1+m. 
\end{equation}
The order $|t_1|=|i|$ is odd by (b), implying that $|t_1| \ne |t_1+m|,$ and so $si = t_1$ holds in \eqref{CI}.
We conclude that $s R = T,$ so $R$ is a CI-subset.  The lemma is proved. 
\QED

\begin{lem}\label{Lem8}
Let $n = 2m$ and $S =\{0,u,v,v+m\}$ such that 
\begin{enumerate}[(a)]
\item $S$ generates $\Z_n;$
\item $1< u < n,$ $ u \mid n$ but $u \nmid m;$ 
\item $\aut(H(\Z_n,S))_{0^+}$ leaves the set $\{0^-,u^-\}$ setwise fixed.
\end{enumerate}
Then $S$ is a BCI-subset.
\end{lem}

\proof 
Let $\delta$ be the partition of $V$ defined in \eqref{D}. 
Applying Lemma \ref{Lemnew} with $R=S$, $T=\{0,u\}$ 
and $t=0$, we obtain that $\delta$  is a system of blocks for $A = \aut(H(\Z_n,S))$. Thus the stabilizer 
$A_{0^+}$ leaves the set $V_0$ setwise fixed, and 
we may consider the action of $A_{0^+}$ on $V_0$. 
The subgraph of $H(\Z_n,S)$ induced  by the set $V_0$ is a circuit of length $2 n/u$, thus  
$A_{0^+}$ fixes also the vertex on this circuit antipodal to $0^+$. 
We find that this antipodal vertex is $(u/2+m)^-$. Therefore, $A_{0^+} = A_{(m+u/2)^-},$ and thus 
$S$ is a BCI-subset if and only if $S-u/2+m$ is a CI-subset of $\Z_n$, see Proposition \ref{Prop1}. 
The latter set is 
\[ S - u/2+m = \big\{ \, u/2+m, -u/2+m, v-u/2, v-u/2+m \, \big\}.\]
Since $u \nmid m$, $u$ is even and the order $|u/2+m|$ is odd. Lemma \ref{Lem7} is applicable to 
the set $S - u/2+m$  (choose $i = u/2+m$ and $j = v-u/2$), it gives us that 
$S-u/2+m$ is a CI-subset. This completes the proof.  \QED 

\bigskip

\noindent{\sc Proof of Theorem \ref{MAIN3}.}  Let $S$ be the subset of $\Z_n$ given in Theorem \ref{MAIN3}.  
We deal first with the case when the canonical bicyclic group $C$ is normal in $A = \aut(H(\Z_n,S))$. 

\medskip

\noindent {\sc Case 1.} $C \trianglelefteq A$. 

\medskip

By Theorem \ref{Thm1} there is a bicyclic group $X$ of $H(\Z_n,S)$ such that $X \ne C$. Since $C \trianglelefteq A,$ 
$X$ is generated by a permutation in the form $c^i \varphi_{r,s,0},$ $r\in \Z_n^*,$ $s \in \Z_n,$ and 
ord$(\varphi_{r,s,0}) \ge 2$. The permutation $\varphi_{r,s,0}$ acts on 
both $\Z_n^+$ and $\Z_n^-$ as an affine transformation. This fact  together with  the connectedness of $H(\Z_n,S)$ 
imply that, $\varphi_{r,s,0}$ acts faithfully on $S^-$.  Thus ord$(\varphi_{r,s,0}) \le 4$. 

Suppose that  ord$(\varphi_{r,s,0}) = 4$. We may assume without loss of generality 
that $S^-$ can be obtained as  $S^- = \{ \, (0^-)^{\varphi_{r,s,0}^j} : j \in \{0,1,2,3\} \, \},$ and so 
$S = \{0,s,(r+1)s,(r^2+r+1)s \}$ and $(r^3+r^2+r+1)s=0$. Since $H(\Z_n,S)$ is connected, $\gcd(s,n)=1,$ and 
$(r+1)(r^2+1)=0$. We find that $(c^i \varphi_{r,s,0})^4$ sends $x^+$ to $(x+r(r+1)(r^2+1)i) ^{\, +} = x^+$. 
Since $X = \big\langle c^i \varphi_{r,s,0} \big\rangle$ is bicyclic,  $n=4,$ and so $H(\Z_n,S) \cong K_{4,4}$. 
This, however, contradicts that $C \trianglelefteq A$.

Now, suppose that ord$(\varphi_{r,s,0}) = 3$. If $A_{0^+}$ is transitive on $S^-,$ then it must be regular 
\cite[Theorem 4.3]{KovKMW}.
This implies that $S^-$ splits into two orbits under  $A_{0^+}$ with length $1$ and $3,$ respectively. 
Let $s \in S$  such that $\{s^-\}$ is an orbit under $A_{0^+}$. Then $A_{0^+} = A_{s^-}$, and by Proposition \ref{Prop1},  
$S-s$ is not a CI-subset of $\Z_n$. However, in this case the graph $\cay(\Z_n,S-s)$  is edge-transitive, 
and thus $S-a$ is a CI-subset (see \cite[page 320]{Li02}),  which is a contradiction.

Finally, suppose that ord$(\varphi_{r,s,0}) = 2$. If $r=1,$ then $2 \mid n$ and $s = m,$ where $n=2m$. 
This implies that $S^-$ is a union of 
two orbits of $C^m$, we may write  $S = \{0,m,s,s+m\}$. 
The graph $H(\Z_n,S)$ is then isomorphic to the lexicographical product $C_n[K_2^c]$ 
of an $n$-circuit $C_n$ with the graph $K_2^c$, see Figure 3.  It is easily seen that then 
$A_{0^+}$ is not faithful on the set $S^-,$ which is a contradiction. 

\begin{figure} 
\centering
\begin{tikzpicture}[scale=1.5]---
\draw (-2,1) node[above] {\small $0$};
\fill (-2,1) circle (1.5pt);
\draw (-1,1) node[above] {\small $0$};
\draw (-1,1) circle (1.5pt);
\draw (0,1) node[above] {\small $s$};
\fill (0,1) circle (1.5pt);
\draw (1,1) node[above] {\small $s$};
\draw (1,1) circle (1.5pt);
\draw (2,1) node[above] {\small $2s$};
\fill (2,1) circle (1.5pt);
\draw (3,1) node[above] {\small $2s$};
\draw (3,1) circle (1.5pt);
\draw (-2,0) node[below] {\small $m$};
\fill (-2,0) circle (1.5pt);
\draw (-1,0) node[below] {\small $m$};
\draw (-1,0) circle (1.5pt);
\draw (0,0) node[below] {\small $s+m$};
\fill (0,0) circle (1.5pt);
\draw (1,0) node[below] {\small $s+m$};
\draw (1,0) circle (1.5pt);
\draw (2,0) node[below] {\small $2s+m$};
\fill (2,0) circle (1.5pt);
\draw (3,0) node[below] {\small $2s+m$};
\draw (3,0) circle (1.5pt);

\foreach \i in {-2,0,2}
{
\draw (\i,0) -- (\i+0.95,0);
\draw (\i,1) -- (\i+0.95,1);
\draw (\i,0) -- (\i+0.95,0.95);
\draw (\i,1) -- (\i+0.95,0.05);
}

\foreach \k in {0,2}
{
\draw (\k,0) -- (\k-0.95,0);
\draw (\k,1) -- (\k-0.95,1);
\draw (\k,1) -- (\k-0.95,0.05);
\draw (\k,0) -- (\k-0.95,0.95);
}
\end{tikzpicture}
\caption{The lexicographical product $C_n[K_2^c]$.}
\end{figure}
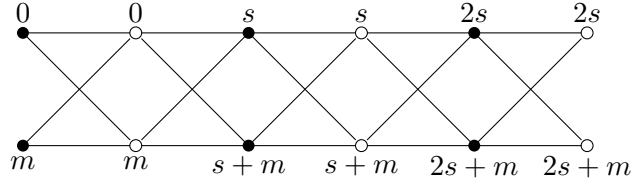

Let $r \ne 1$.  By Lemma \ref{Lem5}, $8 \mid n,$ $r = m+1$ and $s \in \{0,m\},$ where $n=2m$.
We consider here only the case when $s =0$ (the case when $s = m$ can be treated in the same manner). Then $\Z_n^-$ 
splits into the 
following orbits under $\varphi_{r,s,0}$:
\[  \{i^-\}, \; \{ (i+1)^-,(i+1+m)^-\}, \text{ where } i \in \{0,2,\dots,n-2\}.\]
Since $H(\Z_n,S)$ is connected, $S^-$ contains one orbit of length $2$ and two orbits of length $1$. 
Let $S_1$ denote the orbit of length $2$ and let $S_2 = S \setminus S_1$. 
The elements of $S_2$ are even. Let $S_2 = \{s,s'\}$ and let $u = \gcd(s-s',n)$. 
Then $u$ is a divisor of $n$ and also $2 \mid u$. There exist $a \in \Z_n^*$ such that $a (s-s') \equiv u(\mathrm{mod\;} n)$. 
Now, choosing $b = -as'$  (all arithmetic is done in $\Z_n$) ,we find that $a S_2 + b = \{0,u\}$.  
It is clear  that $a S_1 + b = \{v,v+m\}$ for some $v \in \Z_n$.
We finish the proof of this case by showing that the set $R=\{0,u,v,v+m\}$ satisfies the conditions (a)-(c) of Theorem \ref{MAIN2}.  
\medskip

(a): As $H(\Z_n,S)$ is connected, $H(\Z_n,R)$ is also connected.
This implies that  $u,v$ generate $\Z_n$.

\medskip

(c): Since $C \trianglelefteq A,$ $C \trianglelefteq \aut(H(\Z_n,R))$. 
To the contrary assume that $\aut(H(\Z_n,R))_{0^+}$ does not leave $\{0^-,u^-\}$ fixed. Thus there exists  
some $g \in A_{0^+}$  which  maps $v^- $ into $\{0^-,u^-\}$. Letting $w_1^- = (v^-)^g$ and  
$w_2^- = ((v+m)^-)^g,$ we find that $w_1 - w_2 = m,$ and from this that $u=m$. 
However, then  $H(\Z_n,R) \cong C_n[K_2^c],$ which we have already excluded above.
Thus $\aut(H(\Z_n,R))_{0^+}$ fixes setwise $\{0^-,u^-\}$.

\medskip

(b): We have already showed (see previous paragraph) that $u \ne m$ and $1 < u$. 
Since $S$ is not a BCI-subset, $R$ is also a not a BCI-subset. This also implies that $u \mid m$ by Lemma \ref{Lem8}, and we 
conclude that $1 < u < m$ and $u \mid m$, as required. 

\medskip

\noindent {\sc Case 2.} $C \not\trianglelefteq A$. 
	
\medskip

Let $A_{0^+}$ act transitively on $S^-$. This gives that $H(\Z_n,S)$ is edge-transitive. 
Since $C \not\trianglelefteq A,$ $D \not\trianglelefteq A,$ in other words, $H(\Z_n,S)$ is non-normal as a Cayley 
graph over the dihedral group $D$. We apply \cite[Theorem 1.2]{KovKM10}, and obtain that $H(\Z_n,S)$ is either isomorphic to 
$K_n[K_2^c],$ or to one of $5$ graphs of orders $10,14,26,28$ and $30,$ respectively. A direct check by computer shows that 
none of these graphs is possible as all corresponding subsets are BCI-subsets. 

The set $S^-$ cannot split into two orbits under $A_{0^+}$ of size $1$ and $3,$ respectively (see the argument above).
We are left with the case that $S = S_1 \cup S_2,$ $|S_1|=|S_2|,$ and $A_{0^+}$ leaves both sets 
$S_1$ and $S_2$ fixed. For $i \in \{1,2\},$  let $n_i = |\sg{S_i - S_i}|,$ $n_1 \le n_2$, where 
$S_i - S_i = \{a - b\,  : \, a,\, b \in S_i \}$.

We claim that $n_1 = 2$. To the contrary assume that $n_1 > 2$. 
We prove first that $C^{n/n_1} \trianglelefteq A$. 
Applying Lemma \ref{Lemnew} with $R=S$, $T=S_1$ 
and $t=s_1 \in S_1$, we obtain that the partition 
\[ \delta = \big\{ X\cup X^{\psi_{1,s_1,-s_1}} : X \in \orb(C^{n/n_1},V) \big\}, \]
is a system of blocks for $A$. Let us consider the action of $A_{(\delta)}$ (the kernel of $A$ acting on $\delta$) on 
the block of $\delta$ which contains $0^+$.  Denote this block by $\Delta,$ and by $\Delta'$ the block which contains 
$s^-$ for some $s \in S_2$. 
Notice that, the subgraph of $H(\Z_n,S)$ induced by any block of $\delta$ is a circuit of 
length $2n_1,$ and when deleting these circuits, the rest splits into pairwise disjoint circuits of length $2n_2$.   
Let $\Sigma$ denote the unique $(2n_2)$-circuit through $s^-$.
Now, suppose that $g \in A_{(\delta)}$ which fixes $\Delta$ pointwise. 
If $V(\Sigma) \cap \Delta = \{0^+\},$ then $g$ must fix the edge $\{0^+,s^-\},$ and so fixes also $s^-$.  
If  $V(\Sigma) \cap \Delta \ne \{0^+\},$ then $|V(\Sigma) \cap \Delta| = n_2 > 2$. This implies that 
$g$ fixes every vertex on $\Sigma,$ in particular, also $s^-$.   The block $\Delta'$ has at least $n_1$ vertices having a 
neighbor in $\Delta,$ hence by the previous argument we find that all are fixed by $g$. Since $n_1 > 2,$ $\Delta'$ is fixed 
pointwise by $g$. It follows, using the connectedness of 
$H(\Z_n,S),$  that  $g=id_V,$ hence that $A_{(\delta)}$ is faithful on $\Delta$. Thus $C^{n/n_1}$ is a characteristic subgroup of $A_{(\delta)},$ and since $A_{(\delta)} \trianglelefteq A,$ $C^{n/n_1} \trianglelefteq A$. 

Let $G$ be the unique normal subgroup of $A$ that fixes the color classes $\Z_n^+$ and $\Z_n^-$. We consider 
$N = G \cap C_A(C^{n/n_1})$. Then $C \le N$ and $N \trianglelefteq A$.  
Pick  $g \in N_{0^+}$ such that $g$ acts non-trivially on $S^-$. Since $N$ centralizes $C^{n/n_1},$ $g$ fixes pointwise the 
orbit of $0^+$ under $C^{n/n_1},$ and hence also $\Delta$. Then $g^2$ fixes $S^-$ pointwise, and so also $\Delta'$.
We conclude that $g^2 = id_V,$ and that either $N = C,$ or $N = C \rtimes \sg{g}$. 
The case $N = C$ is impossible because $C \not\trianglelefteq A$. 
Let $N = C \rtimes \sg{g}$. Then $(S_i^-)^g=S_i^-$ (for both $i \in \{1,2\}$),  hence 
$S_i$ is a union of orbits of $g$. 
As $g$ normalizes $C,$ $g = \varphi_{r,s,1}$. Recall that ord$(g) =2$. 
If $r \ne 1,$  then by Lemma 3.6, either $C$ is the unique cyclic subgroup of $N,$ or $8 \mid n,$ $r = n/2+1$ and 
$s = 0$ or $s=n/2$. In the former case $C$ is characteristic in $N,$ and since $N \trianglelefteq A,$ $C \trianglelefteq A,$ 
a contradiction. Therefore, we are left with the case that $r=1$ (and so $s=n/2$), or $8 \mid n,$  $r= n/2+1$ and $s=0$ or $s=n/2$.  Then  every orbit of $g$ of length $2$ and contained in $\Z_n^-$ 
 is in the form $\{ j^-,(j+m)^-\}$ as we proved in Case 1. Since $n_i > 2,$ we see that $S_i^-$ must be fixed pointwise by $g$ 
for both $i \in \{1,2\}$. Thus $S^-$ is fixed pointwise by $g,$ implying that $g = id_V$.  This is a contradiction, and so 
$n_1 = 2$.
 
Then $2 \mid n$ and we can write $S_1 = \{s,s+m\},$ where $n = 2m$.
There exist $a \in \Z_n^*$ and $b \in \Z_n$ such that 
$ a S_2 + b = \{0,u\}$ for some divisor $u$ of $n$. Clearly, $a S_1 + b = \{v,v+m\}$ with $v = as + b$. 
We finish the proof of this case by showing that the set 
$R=\{0,u,v,v+m\}$ satisfies the conditions (a)-(c) of Theorem \ref{MAIN2}.  
\medskip

(a): As $H(\Z_n,S)$ is connected, $H(\Z_n,R)$ is also connected.
This implies that  $u,v$ generate $\Z_n$.

\medskip

(c): Since $S_1$ and $S_2$ are left fixed by $A,$ $\aut(H(\Z_n,R))_{0^+}$ leaves the set $\{0^-,u^-\}$ fixed.

\medskip

(b): If $u=1,$ then $\aut(H(\Z_n,\{0,u\})) \le  D_{4n}$. But then $C \trianglelefteq A,$ which is a contradiction. 
We conclude that $1 < u$, and by Lemma \ref{Lem8} that $u \mid m$ also holds, i.e., $1 < u < m$ and $u \mid m,$ as 
required. \QED

\end{document}